\documentclass[12pt]{article}
\usepackage{fancyvrb}

\usepackage[breakable]{tcolorbox}
 \usepackage{parskip} 
    
\usepackage{amsfonts,amssymb,amsthm}
\usepackage{physics}
\usepackage{graphicx}
\usepackage{epstopdf}
\usepackage{algorithmic}
\ifpdf
  \DeclareGraphicsExtensions{.eps,.pdf,.png,.jpg}
\else
  \DeclareGraphicsExtensions{.eps}
\fi

\newtheorem{remark}{Remark}
\newtheorem{proposition}{Proposition}

\newtheorem{example}{Example}
\newtheorem{definition}{Definition}
\newtheorem{theorem}{Theorem}
\newtheorem{lemma}{Lemma}

\usepackage{amsopn}

\usepackage{color,soul}

\usepackage{physics}
\usepackage{lineno}
\usepackage{graphics}
\usepackage{}
\graphicspath{{Figures/}}
\usepackage{subcaption}
\usepackage{float}
\usepackage{pgf,tikz}
\usetikzlibrary{arrows,calc,through,matrix,decorations.pathreplacing,angles,arrows.meta,positioning,shapes,quotes}
\usepackage{hyperref}
\usepackage{enumitem}
\usepackage{pifont}
\usepackage{textcomp} 
\usepackage{listings}
\usepackage{xcolor}
\definecolor{codegreen}{rgb}{0,0.6,0}
\definecolor{codegray}{rgb}{0.5,0.5,0.5}
\definecolor{codepurple}{rgb}{0.58,0,0.82}
\definecolor{backcolour}{rgb}{0.95,0.95,0.92}
\lstdefinestyle{mystyle}{
    backgroundcolor=\color{backcolour},   
    commentstyle=\color{codegreen},
    keywordstyle=\color{magenta},
    numberstyle=\tiny\color{codegray},
    stringstyle=\color{codepurple},
    basicstyle=\ttfamily\footnotesize,
    breakatwhitespace=false,         
    breaklines=true,                 
    captionpos=b,                    
    keepspaces=true,                 
    numbers=left,                    
    numbersep=5pt,                  
    showspaces=false,                
    showstringspaces=false,
    showtabs=false,                  
    tabsize=2
}

\lstdefinelanguage{Maple}%
{morekeywords={and,assuming,break,by,catch,description,do,done,%
elif,else,end,error,export,fi,finally,for,from,global,if,%
implies,in,intersect,local,minus,mod,module,next,not,od,%
option,options,or,proc,quit,read,return,save,stop,subset,then,%
to,try,union,use,uses,while,xor},%
sensitive=true,%
morecomment=[l]\#,%
morestring=[b]",%
morestring=[d]"%
}[keywords,comments,strings]%
\newcommand{\red}{\color{red}{}}
\definecolor{bgreen}{rgb}{0.0, 0.7, 0.0}
\definecolor{Tim}{rgb}{.5,0.15,.02}
\definecolor{jim}{rgb}{0,.4,0}
\definecolor{sana}{rgb}{0.7,0,0.7}
\newcommand{\black}{\color{black}{}}
\newcommand{\blue}{\color{blue}{}}
\definecolor{LightGray}{gray}{0.9}
\newcommand{\jim}{\color{black}}
\newcommand{\Tim}{\color{black}}
\newcommand{\sana}{\color{black}}
\newcommand{\minimax}{minimax\ }
\newcommand{\ie}{\textit{i.e.}}
\newcommand{\BF}{\boldmath}

\newcommand{\R}{{\mathbb{R}}}

\newcommand{\fm}{{\leftarrow}}
\newcommand{\bm}{\boldmath}
\def\cF{{\cal F}}
\def\cL{{\cal L}}
\def\M{\mathsf{K}}

\def\cSx{{\text{SolSet}(x)}}
\def\cSp{{\text{SolSet}(p)}}

\newcommand{\beq}{\begin{linenomath}\begin{equation}\begin{aligned}} 
\newcommand{\eeq}{\end{aligned}\end{equation}\end{linenomath}} 

\definecolor{cmarkgreen}{rgb}{0.0, 0.5, 0.0}
\definecolor{xmarkred}{rgb}{0.8, 0.0, 0.0}
\newcommand{\xmark}{\color{xmarkred}{{\ding{55}}}}
\definecolor{ballblue}{rgb}{0.63, 0.79, 0.95}
\definecolor{bgreen}{rgb}{0.0, 0.7, 0.0}
\tikzset{rdnode/.style={draw,circle,color=red,fill=red!30,text=black}}
\tikzset{bdnode/.style={draw,circle,color=blue,fill=ballblue,text=black}}
\tikzset{bluenode/.style={draw,circle,color=blue,scale=.8,fill=ballblue,text=black}}
\tikzset{rednode/.style={draw,circle,scale=.8pt,color=red,fill=red!30,text=black}}
\tikzset{regnode/.style={draw,circle,scale=.8}}
\tikzset{uninode/.style={draw,circle,scale=.75pt}}
\tikzset{2uninode/.style={draw,circle,scale=.92pt}}
\tikzset{secnode/.style={draw,circle,scale=.75pt,color=red,fill=red!30,text=black}}
\tikzset{2secnode/.style={draw,circle,scale=.92pt,color=red,fill=red!30,text=black}}
\tikzset{tednode/.style={draw,circle,scale=.75pt,color=blue,fill=ballblue,text=black}}
\tikzset{2tednode/.style={draw,circle,scale=.92pt,color=blue,fill=ballblue,text=black}}

\title{Structured Systems of Nonlinear Equations}
 \author{Sana Jahedi\footnote{Department of Mathematics and Statistics, University of New Brunswick. Corresponding author's email: s.jahedi@unb.ca }, Timothy Sauer\footnote{Department of Mathematical Sciences, George Mason University} and James A. Yorke\footnote{IPST, Mathematics, and Physics, University of Maryland College Park.}}

\date{}

\begin{document}

\maketitle

\begin{abstract}
In a ``structured system'' of equations, each equation depends on a specified subset of the variables. 
In this article, we explore properties common to ``almost every'' system with a fixed structure and how the properties can be read from the corresponding connection graph.

A solution $p$ of a system $F(p)=c$ is called robust if it persists despite small changes in $F$. 
We establish methods for determining robustness that depends on the structure, as expressed in the properties of the corresponding directed graph of the structured system. 
The keys to understanding linear and nonlinear structured systems are subsets of variables that we call forward and backward bottlenecks.
In particular, when robustness fails in a structured system, it is due to the existence of a unique ``backward bottleneck’’, that we call a ``minimax bottleneck''.  We present a numerical method for locating the minimax bottleneck. 
We show how to remove it by adding edges to the graph.
\end{abstract}

{\bf keywords}: 
  nonlinear equations, structured systems, robustness, prevalence, generic rank, constant rank theorem, implicit function theorem, dilation

\section{Introduction}
This paper aims at re-framing and extending what is referred to in the engineering literature as ``structured systems'': each equation depends on a specified subset of the unknown variables. {\Tim Substantial previous work on structured systems has been done in the context of}
control theory~\cite{liu_controllability_2011,liu_observability_2013}. We remove it from that framework to broaden its applicability in the sciences.

In mathematical applications, there is  often significant uncertainty in the specific details of the model equations. Yet some basic properties of the model {\Tim are}  independent of these details, at least for generic implementations of the underlying structure of the equations. This fact motivates a focus on structured systems, which form a vector space of systems for each fixed structure. We will show how to guarantee the robustness of solutions for typical systems with a given structure. We find that lack of robustness is equivalent to the absence of a bottleneck in the graph-theoretic representation of the system. In addition, we characterize multidimensional solution sets in the case that solution sets are manifolds.

\begin{definition}\label{def:robust}\normalfont
Let $F:U\subset \mathbb{R}^N\to \mathbb{R}^M$ be a $C^1$ function and assume $U$ is an open subset of $\mathbb{R}^N$. Write $F(x)=(f_1(x),\cdots,f_M(x))$. 
We say a point $x$ is {\bf robust} for $F$
if $DF(x)$ has rank $M$, and $x$ is
{\bf fragile} if it is not robust.  For $p$ in $U$, we define the solution set
\begin{align}\label{def:solset}
  \cSp:=\{x \in U: F(x)=F(p)\}. 
\end{align} 
to be robust if every $x \in \cSp$ is robust. 
\end{definition}

\begin{example}\rm 
As the first example of a robust structured system, consider a model of a biological feedback loop~\cite{murray1989mathematical} in which the messenger mRNA activates an enzyme that interacts with the substrate to create a product that in turn, modulates the expression of mRNA. This negative feedback loop can be modeled by the system
  \begin{eqnarray} \label{eqMEPd}
   \dot{x_1} &=& -c_1+f_1(x_1,x_3) \nonumber\\
 \dot{x_2} &=& -c_2+f_2(x_1,x_2) \\
 \dot{x_3} &=& -c_3+f_3(x_2,x_3), \nonumber
  \end{eqnarray}
  where $x_1, x_2$, and $x_3$ denote levels of mRNA, enzyme, and product, respectively.
  Our interest in this article is in {\it finite dimensional systems of equations}.
  These often arise as equations of steady states of {\it differential} equations.

  The steady states of the above system are solutions of the system of equations $F(x) = c$, or in component form,
   \begin{eqnarray} \label{eqMEP}
c_1 &=& f_1(x_1,x_3) \nonumber\\
c_2 &=& f_2(x_1,x_2) \\
c_3 &=& f_3(x_2,x_3). \nonumber
  \end{eqnarray}
  This system is a structured system;
the form of the equations tells us that the first equation cannot depend upon $x_2$, nor the second upon $x_3$, nor the third upon $x_1$.
Figure~\ref{figRNA}(c) shows a directed graph representing the system~\eqref{eqMEP}, where there is an arrow from node $j$ to node $i$ if variable $j$ is in equation $i$. 
\end{example}

 \begin{figure}[H]
    \centering
     \begin{tikzpicture}[ultra thick]
      \coordinate(o) at (0,0);
     \node at (0,0.5) {$\begin{array}{r}
f_1(x_1,x_3)=c_1,\\
f_2(x_1,x_2)=c_2,\\
f_3(x_2,x_3)=c_3.
\end{array}
  $};
  \node at (4,0.5){$DF = \left[
\begin{array}{cccc}
f_{11}& 0&f_{13}  \\
f_{21}&f_{22}&0\\
0&f_{32} &  f_{33} 
\end{array}
\right]$};
\node[regnode](s3)[xshift=11.6cm,yshift=1.5cm]{$\bm 1$};
\node[regnode](s1)[xshift=9.9cm,yshift=0cm]{$\bm 2$};
\node[regnode](s2)[xshift=13.3cm,yshift=0cm]{$\bm 3$};   
\draw[<-,>=stealth,black](s1)-- (s3);
\draw[->,>=stealth,black](s2)-- (s3);
\draw[->,>=stealth,black](s1)-- (s2);
\draw[->,>=stealth,black](s3)edge[in=-20,out=60,loop above]node[below right]{}();
\draw[->,>=stealth,black](s2)edge[in=20,out=-60,loop right]node[below left]{}();
\draw[->,>=stealth,black](s1)edge[in=-20,out=60,loop left]node[below left]{}();
\path ([xshift=-27ex,yshift=-2ex]current bounding box.south) node[text width=4em] (sca){\subcaption{ }};
 \path node[text width=5em,right=6em of sca] (scb){\subcaption{ }};
\node[text width=5.5em,right=5em of scb] (scc){\subcaption{ }};
    \end{tikzpicture} 
    \caption{ {\sana Three representation types of a feedback loop involving $x_1$, mRNA density, $x_2$, enzyme density, and $x_3$, the product density.} (a) Structured system of three equations in three variables. (b) The Jacobian matrix has rank 3 for generic entries. (c) Directed graph corresponding to the system in (a).}\label{figRNA}
\end{figure}
The property of robustness, that a small variation in the model equations does not destroy the steady state solution, is extremely desirable for models of physical systems.
If a solution exists for a certain choice of parameters, but does not exist for nearby values, it is likely that the model will not match real behavior in nature when parameter uncertainty, as well as model uncertainty, are likely to be significant.  Mathematically speaking, the continued existence of a solution under small uncertainty can often be guaranteed by the Implicit Function Theorem~\cite{rudin1976principles}. For the system~\eqref{eqMEP}, if the Jacobian of the system at the solution is full rank, the Implicit Function Theorem says that the solution extends to small perturbations of the system. This local result can be tested on a case-by-case basis if enough detailed information about the system is known.

A more global approach was taken in the paper~\cite{JahediRobust}, from which follows the following result; 
for almost every function $F$ of the form~\eqref{eqMEP} and almost every $p$ in the domain $\R^3$, the Jacobian $DF(x)$ at a point $x$ in SolSet$(p)$ achieves the maximum possible rank of 3. This ``global'' result shows that the rank of $DF(x)$ is typically 3, and so the Implicit Function Theorem will apply for almost every instantiation of~\eqref{eqMEP}, and allow us to conclude the solutions are robust. Of course, such a global statement hinges on one's definition of ``almost every''. Since these function spaces of $C^\infty$ functions are infinite-dimensional in general, the ``almost every $F$'' in this result is in the sense of prevalence~\cite{sauer_embedology_1991,hunt_prevalence_1992}.
 
More generally, the following facts were proved in~\cite{JahediRobust}. Let $U\subset \R^N$ be an open set, and let $\cF$ be any vector space of $C^\infty$ functions $F:U\subset \R^N\to \R^M$ respecting a structure, meaning that only certain variables are allowed to appear in each equation. The classical case where ``only certain variables'' means ``all variables'' is also allowed. In the above case, that means all functions $F$ in $\cF$ have form~\eqref{eqMEP}. Assume also that $\cF$ contains all linear functions that respect the structure. Due to this assumption and the fact that $U$ is open, there is a maximum possible rank of the Jacobian of the system on $U$, which depends only on the structure; we denote this rank by $r$. (See Definition \ref{def21} for a rigorous definition.) Then Theorem 2.13 of \cite{JahediRobust} says the following.

\begin{theorem}\label{thm:CR} {\rm \cite{JahediRobust}}
For almost every $F\in\cF$ and almost every $p\in U$, the following holds:

 SolSet$(p)$  is a $C^\infty$-manifold of dimension $N-r$ and for all $x\in\cSp$, $DF(x)$ has kernel of dimension $N-r$.
 \end{theorem}
The following two propositions follow immediately from Theorem \ref{thm:CR} and the definitions of robust and fragile. 
\begin{proposition}\label{Prop:key} If $r=M$, then SolSet$(p)$ is robust for almost every $F\in\cF$, and almost every $p\in U$.
 \end{proposition}

 \begin{proposition}\label{prop2}
 If $r<M$, then  $p$ is fragile for every $F$ and every $p\in U$. In particular, SolSet$(p)$ is fragile for every $F$ and every $p\in U$.
  \end{proposition}

Theorem~\ref{thm:CR} and Propositions~\ref{Prop:key} and \ref{prop2} can be considered as a global extension of the Implicit Function Theorem.
The Implicit Function Theorem describes the solutions in a neighborhood of $p$, whereas the above results are global results applicable to almost every $p$. These global results reveal the key importance of the ``maximum possible'' rank of the Jacobian in a structured system. {\sana In the literature, this maximum possible rank of a matrix is often referred to as the generic rank.} In this article we establish a { graph-theoretic criterion} that is equivalent to the rank $DF(x) = M$ for Lebesgue-almost every $x$ and almost every $F:U\subset \R^N \to \R^M$ in the function space of a structured system.

Next we discuss a system that does not have {\sana a robust solution}.

\begin{example} \rm \label{exCEP}
 Consider an ecological system  where the logarithmic growth rates of species 1 and 2 depend only upon species 3,
\begin{eqnarray} \label{CE1}
\frac{\dot{x}_1}{x_1}  &=&-c_1+f_1(x_3), \nonumber\\
\frac{\dot{x}_2}{x_2}  &=&-c_2+f_2(x_3),\\
\frac{\dot{x}_3}{x_3}  &=&-c_3+f_3(x_1,x_2,x_3),\nonumber
\end{eqnarray}
where $x_1$ and $x_2$ denote the population densities of the predators, and $x_3$ is the prey.
 To search for steady states of the system, set the left sides of the equations to zero, yielding the system of equations
  \begin{eqnarray} \label{eqCEP}
 f_1(x_3) &=& c_1,\nonumber\\
 f_2(x_3) &=& c_2,\\
 f_3(x_1,x_2,x_3) &=& c_3.\nonumber
  \end{eqnarray}
See Figure~\ref{fg:Fragile} for other representations of the above structured system of equations. The first two equations share one unknown, $x_3$. 
 There may be a solution of these equations with positive $x_1, x_2, x_3$ for some exceptional $c=(c_1,c_2,c_3)$, but it will ``fail to be robust'',  in the sense that for almost every choice of $f_1$ and $f_2$, there will be no solutions, since the roots of $f_1$ and $f_2$ will fail to overlap. This lack of robustness is an example of a concept long known in the ecological literature as the Competitive Exclusion Principle.
 \end{example}

 The property that two of the equations $f_1, f_2$ in~\eqref{eqCEP} depend only on one variable is the cause of the non-robustness in this example. More generally, when a structured system has a subset of $m$ equations that collectively depend on fewer than $m$ variables, we say there is a ``backward bottleneck'' in the system (see the next section for a more precise definition). 
 The main result of this article is that if the maximum possible rank $r$ is less than $M$, then a structured system must have such a backward bottleneck.
 In fact, this can be viewed as a generalization of the phenomenon expressed by the Competitive Exclusion Principle \cite{hardin1960}.
 Conversely, in the absence of such a bottleneck, Proposition~\ref{Prop:key} holds, and almost every solution is robust. 
 
 { Backward bottlenecks, called ``dilations'' in a slightly different context, were studied for structural observability and controllability by C.T. Lin \cite{lin_structural_1974}, and extended by Liu et al. \cite{liu_controllability_2011,liu_observability_2013}.  We use backward bottlenecks for a different purpose in this work. Furthermore, we introduce the concepts of ``forward bottleneck'' and ``minimax bottleneck''.
 In addition, we show how bottlenecks can be located computationally from the knowledge of the structure alone.}
 
 \begin{figure}[H]
    \centering
     \begin{tikzpicture}[ultra thick]
  \coordinate (o) at (0,0);
     \node at (0,0) {$\begin{array}{r}
f_1(x_3)=c_1,\\
f_2(x_3)=c_2,\\
f_3(x_1,x_2,x_3)=c_3.
\end{array}
  $};
\path ([yshift=-2ex]current bounding box.south) node[text width=4em] (sca){\subcaption{ }};
\node[rednode](s3)[xshift=11.5cm]{\bf 3};
\node[bluenode](s1)[xshift=9.5cm]{\bf 1};
\node[bluenode](s2)[xshift=13.5cm]{\bf 2};    
\draw[<->,>=stealth,black](s1)-- (s3);
\draw[<->,>=stealth,black](s2)-- (s3);
 \path node[text width=5em,right=6.7em of sca] (scb){\subcaption{ }};
\node[text width=5em,right=4.4em of scb] (scc){\subcaption{ }};
\draw[->,>=stealth,black](s3)edge[in=-20,out=60,loop above]node[below right]{}();
\node at (4.3,0){$DF = \left[
\begin{array}{cccc}
0& 0&f_{13}  \\
0 & 0&f_{23}\\
f_{31} & f_{32}& f_{33} 
\end{array}
\right]$};
    \end{tikzpicture} 
    \caption{ \small
    {\bf 
    A fragile  structured system of equations motivated by Competitive Exclusion Principle.} {\bf (a)} A structured system of equations describing the steady states of system~\eqref{eqCEP}. For example, ${f_1}$ in the first equation is allowed to depend on ${x_3}$, but not ${x_1}$ or ${x_2}$. 
    This fact is represented in the two other parts of this figure. {\bf (b)} The structure matrix $DF(x) =\big[\frac{\partial f_i }{\partial x_j}(x)\big]$. {\bf (c)} The directed graph of the system. An edge from node $i$ to node $j$ in the graph means that variable $i$ is allowed to appear in equation $j$.
    Systems of this form cannot have a robust solution, so any solution that exists is fragile. 
    The coloring in (c) illustrates a ``backward bottleneck''. The ``bottle'' nodes 1 and 2 are blue, and the ``neck'' node 3 is red here and throughout the paper.
    }
  \label{fg:Fragile}
\end{figure}

\begin{example}\rm
 The {\sana fragile} solution of the ecological system~\eqref{eqCEP} can be made robust.
Adding another prey species (node 4) to  Figure~\ref{fg:Fragile}(c) yields the directed graph in Figure~\ref{fg:rbstExpl}(c), and the equations
  \begin{eqnarray} \label{eqRob}
 f_1(x_3,x_4) &=& c_1,\nonumber\\
 f_2(x_3,x_4) &=& c_2,\\
 f_3(x_1,x_2,x_3) &=& c_3,\nonumber\\
 f_4(x_1,x_2) &=& c_4.\nonumber
  \end{eqnarray}
Let $\cF$ be the vector space of all $C^1$ functions $F=(f_1,f_2,f_3,f_4)$ where the $f_i$ are restricted to the form shown by the structured system of equations~\eqref{eqRob}. 

We will find that {\sana the graph} in Figure~\ref{fg:rbstExpl}(c) has no backward bottleneck. Let $\cF$ denote the set of all $C^\infty$ functions respecting the structured system of equations~\eqref{eqRob}. For almost every $F\in\cF$ and for almost every $x=(x_1,\ldots,x_4)$ that is a solution of the system in Figure~\ref{fg:rbstExpl}(a),  each sufficiently small perturbation of {\sana $F$} also has a solution. Thus, such solutions are allowed to exist in naturally-occurring circumstances.
\end{example}
\begin{figure}[H]
    \centering
     \begin{tikzpicture}[ultra thick]
  \coordinate (o) at (0,0);
     \node at (0,0) {$\begin{array}{r}
f_1(x_3, x_4)=c_1,\\
f_2(x_3, x_4)=c_2, \\
f_3(x_1, x_2,x_3)=c_3,\\
f_4(x_1, x_2)=c_4.
\end{array}
  $}; 
   \path ([yshift=-2ex]current bounding box.south) node[text width=5em] (sca){\subcaption{ }};
 \node[regnode](s1)at(7.8cm,0.57cm){\bf 1};
 \node[regnode](s2)at(11cm,0.57cm){\bf 2};
\node[regnode](s3)at(9.4cm,0.57cm){\bf 3};
 \node[regnode](s4)at(9.4,-0.5){\bf 4};

\draw[->,>=stealth,black](s3)edge[in=-20,out=60,loop above]node[below right]{}();
\draw[<->,>=stealth,black](8.1cm,0.57)-- (9.1cm,0.57);
\draw[<->,>=stealth,black](9.7cm,0.57)-- (10.66cm,0.57);
\draw[<->,>=stealth,black](8cm,0.35)-- (9.12cm,-0.45);
\draw[<->,>=stealth,black](10.8cm,0.35)-- (9.68cm,-0.45);
  \path node[text width=3em,right=8em of sca] 
  (scb){\subcaption{ }};
\node[text width=3em,right=6.2em of scb]
(scc){\subcaption{ }};
   \node at (4.6,0){$DF = \left[
\begin{array}{cccc}
0& 0 &f_{13} &  \bm{f_{14}}\\
0 & 0 &\bm{f_{23}} & f_{24}\\
\bm{f_{31}}& f_{32} & f_{33} & 0\\
f_{41} &  \bm{f_{42}} & 0 & 0
\end{array}
\right]$};
    \end{tikzpicture} 
    \caption{\small {\bf A robust family of systems}.
   For almost every $F$ {\sana that respects the given structure in (a)}, $DF(x)$ is nonsingular, so the system will have ``robust'' solutions.  (a) The structured system of equations. (b) The Jacobian is generically of rank 4. (c)  No bottleneck exists in the associated directed graph.
}\label{fg:rbstExpl}
\end{figure}

In 1974 C.T. Lin~\cite{lin_structural_1974} developed a theory of linear structured systems for observability and control purposes.  Lin's ideas were further developed in the control and observability literature, often without proofs; see \cite{shields1976structural,blackhall2010structural,liu_controllability_2011,liu_observability_2013,dion_generic_2003,liu_graphical_2019} for examples. In this article, we study a different, and in a way, more fundamental question: when are solutions of a (nonlinear) structured system robust?

We extend Lin's ideas to nonlinear structured systems and rigorously prove the connection between bottlenecks and the maximum rank of the Jacobian at solutions. In addition, we show that the principal bottleneck in a system, called the minimax backward bottleneck, can be located through the concept of ``kernel nodes'', and present a computational approach to identify them.

\section{Structured systems and directed graphs}\label{sec:structuredsystem}
A convenient way to visualize a structured system is to assign a directed graph to the system.
 For a general system of $M$ equations in $N$ variables, we consider a graph of $P$ nodes where $P=\max\{M,N\}$. Such graphs are illustrated in  Figures~\ref{figRNA}(c), \ref{fg:Fragile}(c) and~\ref{fg:rbstExpl}(c)  where $M=N$; several examples later in this section treat cases where $M\neq N$.

{\bf Graph Assumptions.} 
In graph $G$, having an edge from node $j$ to node $i$ means that 
variable $j$ is allowed to appear in function $f_i$. 
We call such a node $i$ (having an incoming edge) a function node and such a node $j$ (having an outgoing edge) a variable node. We assume that each node has either at least one incoming edge or at least one outgoing edge. 
In particular, $i\leq M$ means node $i$ has an incoming edge, and $j \leq N$ means node $j$ has an outgoing edge.
A node $k$ satisfies $k\leq \min{(M,N)}$ if and only if it is both a function and a variable node. Such a labeling can always be achieved by numbering all variables in the system arbitrarily from $1$ to $N$, and then listing the equations in arbitrary order. Note that although every system of equations realizes a directed graph in this way, some directed graphs cannot be realized from a system of equations (such as $1\rightarrow 2$ or $2\rightarrow 1$).

For a directed graph $G$ of a structured system, define the vector space {\BF $\cF(G)$} to be the space of all $C^\infty$ functions $F:U\subset \R^N\to \R^M$ that {\bf respect the graph $G$}, in the sense that the variable $j$ is allowed to appear in equation $i$ only if there is an edge from node $j$ to node $i$. 
The subspace $\bm{\cL(G)}$ denotes the vector space of linear functions that respect the graph $G$. 
Other notable subspaces of $\cF(G)$ include the subspace of all polynomial $F$, or polynomials with some maximum degree; \ie, the components $f_i$ are polynomials with some maximum degree.
 In the latter case, the vector space is finite dimensional.
 
 An alternative, more algebraic way to represent a structured system is by the structure matrix of partial derivatives.
 
{
\begin{definition}\label{def21} \rm
The {\bf structure matrix} $S(G)$ of a directed graph is a matrix where the $ij$  entry is allowed to be nonzero if and only if $G$ has an edge from node $j$ to node $i$, such as  illustrated in Figures~\ref{figRNA}(b),  \ref{fg:Fragile}(b) and~\ref{fg:rbstExpl}(b).
A function $F$ is said to {\bf \BF respect a structure matrix $S$} if
$S_{ij}=0$ implies that $\tfrac{\partial F_i}{\partial x_j}(x) = 0$ for all $x$. In particular, let $\mathcal{L}(S)$ be the set of all linear functions $Ax$
where $A$ is a matrix that respects $S$. {\Tim The maximal rank of all matrices respecting a structure matrix $S$ is called {\bf maxrank}$(S)$.}
We say a vector space $\mathcal{F}$ of $C^1$ functions  that respect a structure is a {\bf structured function space}, provided $\mathcal{F}$ includes $\mathcal{L}(S)$.
\end{definition}}

\begin{example}\rm Figure~\ref{fg:12} is a graphical representation of a more complex ecological model with 26 species or nodes.
The reader may find it daunting to determine from the graph whether such a system allows robust solutions.
We will show how to analyze whether such graphs can have robust { steady states}, and in fact, 
no system with this graph can have any robust solutions. We return to this system in Example~\ref{ex26nodes}.  

\end{example}
 \begin{figure}
\begin{center}
\resizebox{11cm}{5cm}{%
 \begin{tikzpicture}
    [ultra thick]
     \tikzset{every node}=[font=\BF]
    \node[2uninode]at(-5.7cm,2cm)(a11){${1}$};
    \node[2uninode]at(-3.5cm,2cm)(a12){${2}$};
    \node[2uninode]at(0cm,2cm)(a13){${3}$};
    \node[2uninode]at(2.2cm,2cm)(a14){${4}$};
    \node[2uninode]at(5cm,2cm)(a15){${5}$};
    \node[2uninode]at(+7cm,2cm)(a16){${6}$};
    \node[2uninode]at(+9cm,2cm)(a17){${7}$};
    \node[2uninode]at(-6.5cm,.2cm)(a21){${8}$};
    \node[2uninode]at(-4.8cm,.2cm)(a22){${9}$};
    \node[uninode]at(-3cm,.2cm)(a23){${10}$};
    \node[uninode]at(-0.7cm,.2cm)(a24){${11}$};
    \node[uninode]at(1.5cm,.2cm)(a25){${12}$};
    \node[uninode]at(3.5cm,.2cm)(a26){${13}$};
    \node[uninode]at(5.5cm,.2cm)(a27){${14}$};
    \node[uninode]at(+8cm,.2cm)(a28){${15}$};
     \node[uninode]at(-5cm,-1.5cm)(a31){${16}$};
    \node[uninode]at(-1.5cm,-1.5cm)(a32){${17}$};
    \node[uninode]at(3cm,-1.5cm)(a33){${18}$};
    \node[uninode]at(6.5cm,-1.5cm)(a34){${19}$};
    \node[uninode]at(-6cm,-3.2cm)(a41){${20}$};
    \node[uninode]at(-4cm,-3.2cm)(a42){${21}$};
    \node[uninode]at(-1cm,-3.2cm)(a43){${22}$};
    \node[uninode]at(1cm,-3.2cm)(a44){${23}$};
    \node[uninode]at(+4cm,-3.2cm)(a45){${24}$};
    \node[uninode]at(6cm,-3.2cm)(a46){${25}$};
    \node[uninode]at(8cm,-3.2cm)(a47){${26}$};
  \draw[<->,>=stealth,black] (a11)-- (a21);
    \draw[<->,>=stealth,black] (a12)-- (a22);
    \draw[<->,>=stealth,black] (a13)-- (a24);
    \draw[<->,>=stealth,black] (a14)-- (a25);
    \draw[<->,>=stealth,black] (a14)-- (a27);
    \draw[<->,>=stealth,black] (a14)-- (a33);
    \draw[<->,>=stealth,black] (a14)-- (a33);
    \draw[<->,>=stealth,black] (a15)-- (a26);
    \draw[<->,>=stealth,black] (a16)-- (a28);
    \draw[<->,>=stealth,black] (a17)-- (a28);
    \draw[<->,>=stealth,black] (a21)-- (a31);
    \draw[<->,>=stealth,black] (a21)-- (a32);
     \draw[<->,>=stealth,black] (a22)-- (a32);
    \draw[<->,>=stealth,black] (a23)-- (a31);
     \draw[<->,>=stealth,black] (a23)-- (a32);
     \draw[<->,>=stealth,black] (a23)-- (a33);
     \draw[<->,>=stealth,black] (a23)-- (a34); 
    \draw[<->,>=stealth,black] (a24)-- (a32);
    \draw[<->,>=stealth,black] (a25)-- (a32);
    \draw[<->,>=stealth,black] (a26)-- (a33);
    \draw[<->,>=stealth,black] (a27)-- (a32);
 \draw[<->,>=stealth,black] (a27)-- (a34);
  \draw[<->,>=stealth,black] (a28)-- (a34);
  \draw[<->,>=stealth,black] (a34)-- (a47);
  \draw[<->,>=stealth,black] (a34)-- (a46);
  \draw[<->,>=stealth,black] (a34)-- (a45);
  \draw[<->,>=stealth,black] (a34)-- (a44);
  \draw[<->,>=stealth,black] (a33)-- (a44);
  \draw[<->,>=stealth,black] (a31)-- (a44);
  \draw[<->,>=stealth,black] (a31)-- (a42);
  \draw[<->,>=stealth,black] (a32)-- (a42);
  \draw[<->,>=stealth,black] (a32)-- (a41);
  \draw[<->,>=stealth,black] (a32)-- (a43);
  \draw[<->,>=stealth,black] (a32)-- (a44);
\end{tikzpicture}
}
\caption{\small A network model adapted from Sol\'e and Montoya~\cite{sole_complexity_2001}. 
Solutions of the structure system of 26 equations in 26 unknowns associated with this graph cannot be robust since a bottleneck exists with nodes 6 and 7. Several other obstructions to robustness exist. See
Figure~\ref{fg:26} for a more complete discussion.}
\label{fg:12}
\end{center}
\end{figure}

{\Tim
\begin{remark}
 According to Proposition \ref{prop2}, if maxrank$(S(G)) = r < M$, then for every function $F$ in the structured function space of $S(G)$, every solution set is fragile. In such a case, we call the graph $G$, or equivalently, the structured matrix $S(G)$, fragile. The main result of this section is Theorem~\ref{thm:BotNck}, which states that the maxrank condition for $S(G)$ to be fragile is equivalent to the existence of a  graphical obstruction that we call a backward bottleneck in $G$.
\end{remark}
}

 Let $G$ be a graph representing a structured system of $M$ equations in $N$ variables.  Let $B$ be a subset of the nodes of $G$. The {\bf forward set of \BF$B$}, denoted by $B^\to$,
is the set of all nodes $g$ in $G$ for which there is an edge starting at a node in $B$ and ending at $g$. The {\bf backward set of \BF$B$}, denoted by $B^\fm$,
is the set of all nodes $g$ in $G$ for which there is an edge starting at $g$ and ending at a node in $B$.

{Let $G_x \subset G$ be the subset of all variable nodes, and let $G_f \subset G$ be the subset of all function nodes.} Let $B$ be a subset of $G_x$, and let $B^\to$ denote its forward set.
For $K >0$, we say the pair of sets  $B$ and $B^\to$ 
 is a {\bf forward ${\bm K}$-bottleneck} if
$B$ has {\bf exactly ${\bm K}$ nodes more} than $B^\to$.
Analogously, we say a pair of sets  of nodes,  $B\subset G_f$ and $B^\fm$, 
 is a {\bf backward $\bm{K}$-bottleneck} if
$B$ has {\bf exactly $\bm{K}$ nodes more} than $B^\fm$.  If $N>M$, there must be a forward bottleneck, and if $M>N$, there must be a backward bottleneck.

The backward bottleneck in Figure~\ref{fg:Fragile}(c) is the pair of sets $B=\{1,2\}$, shaded in blue,  and $B^\fm =B^\to = \{3\}$, shaded in red.   
In Figure~\ref{trophic2}, $B$ is the top level of three nodes, and  $B^\fm=B^\to$ is the lower level of two nodes. In both examples, $B$ is part of both a forward and backward bottleneck.
We refer to $B$ as the {\bf bottle} (and usually color its nodes blue) and $B^\to$  (or $B^\fm$)
as the {\bf neck} (usually colored red).  
It is possible for a node to be both in a bottle and in a neck, in which case we color it both red and blue, as node 1 in examples~\ref{ex:simplefw}~and~\ref{ex:simplebw}.

A system that has a bottleneck will sometimes have many.
 Let $K_\text{max}$ be the largest value $K$ for which there is a forward $K$-bottleneck.  Let $(B, B^\to)$ be a forward $K_\text{max}$-bottleneck where $B$ has as few nodes as possible; i.e., it has the minimum number of nodes that has the maximum $K$. We call such a bottleneck a {\bf minimax} forward bottleneck.   The backward minimax bottleneck is defined analogously.
 It turns out that \minimax backward bottleneck and  \minimax forward bottleneck are unique if they exist (see Theorem~\ref{thm:BotNck} and Remark~\ref{Re:fW-botneck}).
 
The Bottleneck Theorem, Theorem~\ref{thm:BotNck} and Remark~\ref{Re:fW-botneck} below,
provides more detail on bottlenecks that exist in a graph, and shows how that helps detect a bottleneck and how it might be eliminated.
In the proof of the Bottleneck Theorem,
the bottle $B$ we construct is the \minimax bottle.

The set of nodes comprising the minimax bottleneck is not always obvious from the graph. The concept of kernel nodes, described next, allows us to locate the minimax forward and backward bottlenecks. 

\begin{definition} \label{defnn} 
For a matrix $A$, we say a vector $x$ is a {\bf null vector} or a {\bf kernel vector} if $Ax=0$. The kernel of $A$, denoted $\ker A$, is the set of all null vectors of $A$.
For a graph $G$ with structure matrix $S(G)$, the graph nodes (coordinates) can be divided into two distinct types:

(1) {\bf a regular node (or coordinate)} that takes the value zero for every vector in $\ker(A)$ for almost every matrix $A\in \cL(S(G))$, and 

(2) {\bf a kernel node (or coordinate)} that is  nonzero in some null vector of $A$ for almost every $A\in \cL(S(G))$. 

The set of {\bf kernel nodes} is denoted by $B_{\rm kernel}$.
\end{definition}


Kernel nodes can be thought of as  nodes corresponding to coordinates that are nonzero in at least one null vector for a matrix $A$, which is obtained by replacing nonzero entries of the structure matrix $S(G)$ by random numbers.  For example, the matrix in Figure~\ref{fg:Fragile}(b) 
has a null vector $(f_{32},-f_{31},0)^T$
and therefore coordinates 1 and 2 are nonzero in a kernel vector for almost every choice of $f_{ij}$. Hence, nodes 1 and 2 in Figure~\ref{fg:Fragile}(b) are kernel nodes.

Later in this section, we show how to compute kernel nodes by symbolic algebra. Alternatively, one could insert random numbers into the entries allowed by the structured matrix, and compute the kernel. With probability one, the coordinates that are nonzero in kernel vectors will correspond to the kernel nodes. See Appendix B for a worked-out example. For larger problems, this approach may succeed when the symbolic approach becomes excessively computationally complex.

Lemma~\ref{lm20} of Appendix~A implies that each node is either a regular node or a kernel node. To be precise, for each node $j$ of the graph $G$, let $S_j$ denote the set of matrices $A$ in $\cL(G)$ such that the coordinate $x_j = 0$ for all vectors in ker$(A)$. 
According to Lemma~\ref{lm20},  either $S_j$ or its complement is measure zero. Kernel nodes are the ones where $S_j$ has zero measure.

Next, we state the main theorem of this article. The proof is contained in Appendix~A.  The theorem holds that a structured system of equations is fragile if and only if it has a backward bottleneck. In addition, it identifies the bottle of the minimax backward bottleneck as the set of the kernel nodes $B_{kernel}$ of the transpose of the structure matrix.

 \begin{theorem}[Bottleneck Theorem]\label{thm:BotNck}
Let $S=S(G)$ be the $M\times N$ structure matrix of a directed graph $G$.
Let $T=S^T$ and $r=\text{maxrank}(S)=\text{maxrank}(T)$. 

\begin{enumerate}[label=(\Roman*),ref=(\Roman*),leftmargin=*]
\item   $r<M$ if and only if there is a backward $K$-bottleneck for $K>0$. If so, there exists a unique minimax backward $K^*$-bottleneck where $K^* = M-r.$ \label{thm:BotNck-I}
\item The bottle of the \minimax backward bottleneck of $G$ is the set of kernel nodes $B_{kernel}$ of $T$. 
\item {\sana A graph is fragile if and only if there is
a backward bottleneck.}
\end{enumerate}
\end{theorem}

It follows from the proof (see Appendix) that there is a backward $K$-bottleneck if and only if $0<K \leq K^*$ when $K^*$ is positive.

\begin{remark}
\label{Re:fW-botneck} There are analogous results for forward bottlenecks. The proofs require only trivial changes. Let $S=S(G)$ be the $M\times N$ structure matrix of a directed graph $G$. Let $r=\text{maxrank}(S)$.

(I) $r < N$ if and only if there is a forward $K$-bottleneck for $K>0$. If so, there exists a unique minimax forward $K^*$-bottleneck where $K^* = N-r$.

(II) The bottle of the minimax forward bottleneck of $G$ is the set of kernel nodes $B_{\rm kernel}$ of structure matrix $S$. 

(III) If there is a forward bottleneck, then according to Theorem 1.3, for almost every $C^\infty$ function $F:U\subset {\mathbb R}^N\rightarrow {\mathbb R}^M$ with directed graph $G$ and almost every $p\in U$, SolSet$(p)$ is a $C^\infty$-manifold of dimension $K^*$.
\end{remark}

\begin{remark}
Theorem~\ref{thm:BotNck} in combination with Propositions~\ref{Prop:key} and \ref{prop2}
 shows that almost every function in a structured function space has one of two possibilities: Either (1) almost every solution is robust, if no backward bottleneck exists, or (2) all solutions are fragile, if a backward bottleneck exists.
\end{remark}

 \begin{example}\rm [The simplest forward bottleneck with $N>M$] \label{ex:simplefw}
{Let $M=1, N=2$, meaning that we have a structured system of one equation in two variables with the form $F(x)= 
f_1(x_1,x_2) = c_1$.} The structure matrix and graph are
\begin{figure}[H]
    \centering
 \begin{tikzpicture}
[ultra thick]
  \node[bluenode](a1)at(2.7cm,0){\bf 2};
    \draw[-,line width=0.15em,color=blue,fill=ballblue] (0.835,0) arc (0:180:0.3);
\draw[-,line width=0.15em,color=red,fill=red!30] (0.23,0) arc (180:360:0.3);
  \node[](a2)at(0.523cm,0){\bf 1};
 \draw[black,->,>=stealth](2.35,0)--(0.85,0);
\draw[->,>=stealth,black](a2)edge[in=-20,out=60,loop above]node[below right]{}();
\node at (-3,0){$S=\left[\begin{array}{cc}
f_{11} &f_{12}
\end{array}
\right]$};
    \end{tikzpicture} 
\end{figure}
Here maxrank$\ = 1, K^*=N-$ maxrank$(S) = 1>0$.
There is a forward bottleneck $\{1,2\}\rightarrow \{1\}$ (shown) but no backward bottleneck. 
The structure matrix has a  null vector  $\nu=(
f_{12},-f_{11})^T$. Since 
both coordinates are symbolically non-zero, they are both kernel nodes and correspond to the coordinates of the minimax forward bottle; See Remark~\ref{Re:fW-botneck}(II). The ``bottle'' nodes are blue and the ``neck'' nodes are red throughout the paper.
\end{example}
\begin{example}[The simplest backward bottleneck with $M>N$] \rm \label{ex:simplebw}
Let $M=2, N=1$. The map $F(x) =(f_1(x_1), f_2(x_1))$
 has 
the structure matrix and graph 

\begin{figure}[H]
    \centering
 \begin{tikzpicture}
[ultra thick]
  \node[bluenode](a1)at(2.7cm,0){\bf 2};
    \draw[-,line width=0.15em,color=blue,fill=ballblue] (0.835,0) arc (0:180:0.3);
\draw[-,line width=0.15em,color=red,fill=red!30] (0.23,0) arc (180:360:0.3);
  \node[](a2)at(0.523cm,0){\bf 1};
 \draw[black,->,>=stealth](0.85,0)--(2.35,0);
\draw[->,>=stealth,black](a2)edge[in=-20,out=60,loop above]node[below right]{}();
\node at (-3,0){$S=\left[\begin{array}{c}
f_{11} \\ f_{21} 
\end{array}
\right]$};
    \end{tikzpicture} 
\end{figure}
Here maxrank$\ = 1, K^*=M-$ maxrank($S) = 1$. By Theorem~\ref{thm:BotNck}, there is a backward bottleneck, which is $\{1\}\rightarrow \{1,2\}$, shown in the graph. 
The transpose of the structure matrix has a null vector $\nu = (f_{21},-f_{11})^T$, i.e. $S^T\nu = 0$.
The nonzero components of $\nu$, and therefore the kernel nodes, are \{1,2\}. This set is the bottle of the minimax backward bottle, as is guaranteed by Theorem~\ref{thm:BotNck}(II). Here, $\nu$ has only 2 coordinates.
Since $N-$ maxrank($S) = 0$, by Remark~\ref{Re:fW-botneck} there is no forward bottleneck.
\end{example}

{\bf A procedure for making a graph robust.}  
Let $K^*= M - \text{maxrank}(S(G)).$ If $K^*>0$, then find the minimax backward bottleneck. We can then reduce $K^*$ by one by adding an edge from any node that is not in the neck to any bottle node. Such edges always exist: For example, each node in the bottle $B$ that is not a neck node can be given a self-edge if necessary. When this process is repeated $K^*$ times, the generic rank becomes equal to $M$.
Note that no nodes are being added to the graph in this procedure, so no new bottlenecks are created. 

{\sana {\bf\BF Computational method for finding a bottle.} 
When dealing with a small structured systems such as the one represented in Figures~\eqref{fg:Fragile}~and~\eqref{trophic2} one could use symbolic algebra to find the bottle of the minimax bottleneck. Computer software packages for symbolic algebra give us a computational means of determining the kernel node set $B_{\rm kernel}$ of the matrix $T$, the transpose of the structure matrix $S$.  
See Eq.~\eqref{symbolic} below to see typical output from Python and Maple using symbolic algebra.
The package computes a basis for the kernel space. When symbolic entries are used in the input matrix, for each kernel node, the software returns non-zero symbolic formulas for at least one of the basis vectors, and it always returns zero for the nodes which are not kernel nodes. For an example, consider the network given in Figure~\ref{trophic2}.

\begin{figure}[H]
\begin{center}
 \begin{tikzpicture}
[ultra thick]
  \tikzset{every node}=[font=\bf]
  \node[bluenode](c1)[xshift=-2cm]{1};
  \node[bluenode](c2)[]{2};
\node[bluenode](c3)[xshift=+2cm]{3};
\node[rednode](r1)[xshift=-1cm,yshift=-2cm]{4};
  \node[rednode](r2)[xshift=+1cm,yshift=-2cm]{5};
\draw[<->,>=stealth,black] (r1)-- (c1);
\draw[<->,>=stealth,black] (r1)-- (c3);
\draw[<->,>=stealth,black] (r1)-- (c2);
\draw[<->,>=stealth,black] (r2)-- (c1);
\draw[<->,>=stealth,black] (r2)-- (c3);
\draw[<->,>=stealth,black] (r2)-- (c2);
\draw[->,>=stealth,black](r1)edge[in=-20,out=60,loop below]node[below right]{}();
\draw[->,>=stealth,black](r2)edge[in=-20,out=60,loop below]node[below right]{}();
\end{tikzpicture}
\end{center}
\caption{{\sana Generalization of the Competitive Exclusion Principle for two trophic levels. This structure is fragile since there exists a backward bottleneck; the pair of sets $B=\{1,2,3\}$ and $B^\fm=\{4,5\}$ form the minimax bottleneck.}
}\label{trophic2}
\end{figure}
\subsection*{In Python}\
\begin{lstlisting}[frame = single, basicstyle=\footnotesize, language=Python]

from sympy import MatrixSymbol, Matrix

f = MatrixSymbol('f', 6, 6) # creat a symbolic matrix to use the symbols later

S = Matrix(5,5,[0,0,0,f[1,4],f[1,5],\
                0,0,0,f[2,4],f[2,5],\  
                0,0,0,f[3,4],f[3,5],\
                f[4,1],f[4,2],f[4,3],f[4,4],0,\
                f[5,1],f[5,2],f[5,3],0,f[5,5]])

T = S.transpose() # calculate the transpose of the structure matrix

B = T.nullspace() # calculate a basis for the nullspace of T

Matrix(B) # print the basis for the null space
\end{lstlisting}
\newpage
\subsection*{In Maple}\
\begin{lstlisting}[frame = single, basicstyle=\footnotesize]{Maple}
# Call the linear algebra package
  with(LinearAlgebra): 
# Create the symbolic matrix S
   S := Matrix(5,5, S=[0, 0, 0, f14, f15, 0, 0, 0, f24, f25,
   0, 0, 0, f34,f35,
   f41, f42, f43, f44, 0,
   f51, f52,  f53, 0, f55]): 
# Calculate the basis for the kernel of the transpose of S
  NullSpace(Transpose(S))
\end{lstlisting}

\begin{equation}
\label{symbolic}
\text{output:}
\left\{\left[\begin{array}{c}
f24 f35-f25 f34\vspace{.1cm}\\
f15 f34-f35 f14\vspace{.1cm}\\f25 f14-f15 f24\vspace{.1cm}\\ 0\vspace{.1cm}\\0\end{array}\right]\right\}
\end{equation}
We have reformatted this output for the readers' convenience. 
The output is a single vector means that the kernel is one-dimensional. Therefore, there exists a  one-bottleneck. The nodes 1, 2, and 3 are the bottle nodes since components 1, 2, and 3 of the null vector {\sana are almost always nonzero.} By looking at the incoming edges of these three nodes, we find that the bottleneck is the pair of sets ($B, B^\fm$) where $B=\{1,2,3\}$ and  $B^\fm = \{4, 5\}$.
To eliminate this bottleneck, it suffices to add an edge connecting one of the nodes of the set $\{1,2,3\}$ to one of the nodes from the same set. 
 As a result, by Theorem~\ref{thm:BotNck}, 
the new structure matrix (corresponding to the system after an edge is added) has rank $r=5$. Hence, 
the resulting network will have robust solutions, according to Proposition~\ref{Prop:key}.

Symbolic computation works well for reasonably small problems, but symbolic algebra scales poorly in the size of the graph. For larger problems, one could replace the nonzero entries of the structure matrix with random numbers and compute a basis for the null space of the structure matrix. An example is provided in Appendix B.

\section*{Examples with $\bm{M=N}$}
 \begin{example}\rm \label{ex:3x3}
We interpret Example \ref{exCEP} in terms of bottlenecks.
Consider the graph and system in Figure~\ref{fg:Fragile}, where $M=N=3$. Here maxrank $=2$ for generic entries.
By Proposition~\ref{prop2}, for every $F$ and $c$, the system
$$
\begin{array}{rcl}
f_1(x_3)&=& c_1\\[0.1cm]
f_2(x_3) &=& c_2\\[0.1cm]
f_3(x_1,x_2,x_3) &=& c_3\end{array}
$$
has no robust solutions, because 
maxrank$(F)$ is generically 2 and is never greater than 2.
If $F$ is $C^\infty$, there will be solutions for at most a measure zero set of $c$, and those are fragile. By {\sana Theorem~\ref{thm:CR}} for almost every function with this structure and for almost every $x\in \R^3$, the set
$\cSx$ is a one-dimensional manifold. The structure matrix corresponding to the above system is
$$
DF(x)= \left[
\begin{array}{rrr}
0 & 0 & \frac{\partial f_1}{\partial x_3} \\
0 & 0  & \frac{\partial f_2}{\partial x_3}\\
\frac{\partial f_3}{\partial x_1} &\frac{\partial f_3}{\partial x_2} & \frac{\partial f_3}{\partial x_3}
\end{array}
\right]
$$
The kernel of $DF(x)$ is one-dimensional and $(\frac{\partial f_3}{\partial x_2},-\frac{\partial f_3}{\partial x_1},0)^T$ is a kernel vector,
illustrating that the bottle of minimax forward bottleneck consists of the kernel nodes, which correspond to coordinates 1 and 2. Hence, the minimax forward bottleneck is $\{1,2\} \to \{3\}$, and it is a $K^*$ minimax bottleneck, where $\M^* = N- \text{maxrank}(S) = 3-2 = 1$.   Note that in this example, the minimax forward bottleneck is also the minimax backward bottleneck. Having $M=N$ does not necessarily imply that every forward bottleneck is also a backward bottleneck; see the following example for an instance. In a graph where all the edges are bi-directional, every forward bottleneck is also a backward bottleneck.
\end{example}
\begin{example}\rm \label{ex:m=n,fisnotsameb}
 Let $M=N=3.$ Consider the system
 $$\begin{array}{rcl}
f_1(x_1)&=& c_1,\\[0.1cm]
f_2(x_1) &=& c_2,\\[0.1cm]
f_3(x_2,x_3) &=& c_3.\\
\end{array}$$
The structure matrix and graph (displaying the (minimax) bottleneck) are:

 \begin{centering}
\begin{tikzpicture}
[ultra thick]
\node[]at(-3.5,1){};
\node[]at(1.1,-1.3){(b) minimax};
\node[]at(1.1,-1.8){forward bottleneck};
 \node[bluenode](a1)at(1.2cm,-0.2){\bf 2};
  \node[regnode](a2)at(-0.5cm,-0.2){\bf 1};
  \draw[->,>=stealth,black](a2)edge[in=-20,out=60,loop above]node[below right]{}();
\draw[-,line width=0.15em,color=blue,fill=ballblue]
  (2.9,-0.2) arc (0:180:0.25);
\draw[-,line width=0.15em,color=red,fill=red!30] (2.4,-0.2) arc (180:360:0.25);
 \node[](a3)[xshift=2.64cm,yshift=-0.2cm]{\bf 3}; 
 \draw[->,>=stealth,black](a3)edge[in=-20,out=60,loop above]node[below right]{}();
  \draw[black,->,>=stealth](-0.15,-0.2)--(0.9,-0.2);
\draw[black,->,>=stealth](1.5,-0.2)--(2.35,-0.2);
\node[]at(5.8,-1.3){(c) minimax};
\node[]at(5.5,-1.8){backward bottleneck};
\draw[-,line width=0.15em,color=blue,fill=ballblue]
  (4.25,-0.2) arc (0:180:0.25);
\draw[-,line width=0.15em,color=red,fill=red!30] (3.75,-0.2) arc (180:360:0.25);
 \node[](a11)at(4cm,-0.2){\bf 1};
\draw[->,>=stealth,black](a11)edge[in=-20,out=60,loop above]node[below right]{}();
\node[bluenode](a22)at(5.7cm,-0.2){\bf 2};
\node[regnode](a33)[xshift=9cm,yshift=-0.2cm]{\bf 3};
 \draw[->,>=stealth,black](a33)edge[in=-20,out=60,loop above]node[below right]{}();
\draw[black,->,>=stealth](4.31,-0.2)--(5.37,-0.2);
\draw[black,->,>=stealth](6,-0.2)--(6.9,-0.2);
\node[]at(-3.5,-1.3){(a) structure matrix};
\node at (-3.5,-0.1){
$S= \left[
\begin{array}{rrr}
f_{11} & 0 & 0 \\
f_{21} & 0 & 0\\
0 & f_{32} & f_{33}
\end{array}
\right]$};
\node[]at(-3.5,-2.2){};
\end{tikzpicture}
 \end{centering}
 
 The minimax forward bottleneck is the pair $(B,B^\to)$ where $B=\{2,3\}$ and $B^\to=\{3\}$, and the minimax backward bottleneck is the pair $(C,C^\fm)$ where $C=\{1,2\}$ and $C^\fm=\{1\}$.
\end{example}

Note that a lesson to be learned from Examples~\ref{ex:3x3}~and~\ref{ex:m=n,fisnotsameb} is that when $M=N$, the minimax forward bottleneck is not necessarily the same as the minimax backward bottleneck. In fact, even when $M=N$, a forward bottleneck may not be a backward bottleneck. {\sana A graph is robust when no backward bottleneck exists}.

\section*{\bf Examples with $M<N$}
\begin{example} \rm 
Let  $M=2$ and $N=3$. Consider the system
$$\begin{array}{rcl}
f_1(x_1)&=& c_1,\\[0.1cm]
f_2(x_1,x_2,x_3) &=& c_2,
\end{array}$$
which has the structure matrix and graph (displaying its (minimax) forward bottleneck):
\begin{center}
\begin{tikzpicture}
[ultra thick]
\node at (-3,0){
$S= \left[
\begin{array}{rrr}
f_{11} & 0 & 0 \\
f_{21} & f_{22} & f_{23}
\end{array}
\right]$};
  \draw[-,line width=0.15em,color=blue,fill=ballblue] 
  (2.5,0) arc (0:180:0.25);
\draw[-,line width=0.15em,color=red,fill=red!30] (2,0) arc (180:360:0.25);
  \node[](a1)at(2.25cm,0){\bf 2};
 \draw[->,>=stealth,black](a1)edge[in=-20,out=60,loop above]node[below right]{}();
  \node[regnode](a2)at(0.523cm,0){\bf 1};
  \draw[->,>=stealth,black](a2)edge[in=-20,out=60,loop above]node[below right]{}();
\node[bluenode](a3)[xshift=5cm,yshift=0cm]{\bf 3}; 
 \draw[black,->,>=stealth](0.85,0)--(1.9,0);
\draw[black,->,>=stealth](3.7,0)--(2.6,0);
\node[]at(-3.5,-0.8){};
\end{tikzpicture}
\end{center}
The graph has no backward bottleneck. In fact, $M=$  maxrank$(\cF)=2$, hence,
the graph is robust. On the other hand, the kernel is one dimensional with a basis vector
$\nu = (0,f_{23},-f_{22})^T,$ where $S\nu=0$. There exists a 1-forward bottleneck. According to Remark  \ref{Re:fW-botneck}, for almost every $C^\infty$ function $F$ that respects the above structure, and for almost every $p$, SolSet($p$) is a one-dimensional manifold.
 
\end{example}
\begin{example}\rm 
Let $M=3$ and  $N=5$, and consider the system 
$$\begin{array}{rcl}
f_1(x_1,x_2, x_3) &=& c_1,\\[0.1cm]
f_2(x_4) &=& c_2,\\[0.1cm]
f_3(x_5) &=& c_3.
\end{array}$$
 The minimax {forward} bottleneck is a 2-bottleneck $\{1,2,3\} \to \{1\}$.{ \sana Hence, by Theorem~\ref{thm:CR}} for almost every $C^\infty$ function $F$ that respects the above structure and for almost every $p$, SolSet($p$) is a two-dimensional manifold.
 Since maxrank $= M=3$, the graph is robust.
\begin{center}
\begin{tikzpicture}[ultra thick]
   \tikzset{every node}=[font=\bf]
  \node[bluenode](c2)at(-0.3,0){2};
\draw[-,line width=0.15em,color=blue,fill=ballblue] 
  (-1.7,-0.7) arc (0:180:0.3);
\draw[-,line width=0.15em,color=red,fill=red!30] (-2.3,-0.7) arc (180:360:0.3);
  \node[](c1)[xshift=-2cm,yshift=-0.7cm]{1};
\node[bluenode](r1)at(-0.3,-1.6){3};
\node[regnode](c3)[xshift=+2cm]{4};
  \node[regnode](r2)[xshift=+2cm,yshift=-2cm]{5};
\draw[->,>=stealth,black] (r1)-- (c1);
\draw[<-,>=stealth,black] (c1)-- (c2);
\draw[->,>=stealth,black] (c3)-- (c2);
\draw[<-,>=stealth,black] (r1)-- (r2);
\draw[->,>=stealth,black](c1)edge[in=-20,out=60,loop left]node[below right]{}();
\end{tikzpicture}
\end{center}
\end{example}
\subsection*{\bf An example with $\bm{M>N}$. See also Example~\ref{ex:simplebw}}

\begin{example}\rm  Let M=3 and N=2.
Consider the following structured system of equations.
$$\begin{array}{rcl}
f_1(x_1,x_2)&=&c_1,\\[0.1cm]
f_2(x_1,x_2)&=&c_2,\\[0.1cm]
 f_3(x_1,x_2)&=&c_3.
 \end{array}$$
The system has no robust solutions -- due to the existence  of a minimax backward bottleneck.
Let $S^T$ denote the transpose of the structure matrix associated with this structured system. 

 $$  S^T=\left[\begin{array}{rcl}
       f_{11} &f_{21}&f_{31}  \\[0.1cm]
       f_{12} &f_{22}&f_{32} 
   \end{array}\right]$$
For almost all choices of $F=(f_1,f_2,f_3)$, the kernel  of $(S^T)$ is one dimensional, and $\nu=\left[\begin{array}{c}
  f_{21} f_{32} - f_{31} f_{22} \\
  f_{31} f_{12}- f_{11}  f_{32} \\
f_{11} f_{22} -  f_{21} f_{12}
 \end{array}\right]$
  is a null vector so $S^T\nu=0$. Therefore, the backward bottle 
  is  $B=\{1,2,3\}$, since all the components of the null vector are nonzero. Hence  $B^\fm$ equals $\{1,2\}$. Therefore the pair $(B,B^\fm)$ is a minimax 1-backward bottleneck:
  \begin{figure}[H]
    \centering
\begin{tikzpicture}
[ultra thick]
  \draw[-,line width=0.15em,color=blue,fill=ballblue]
  (3,0) arc (0:180:0.3);
\draw[-,line width=0.15em,color=red,fill=red!30] (2.4,0) arc (180:360:0.3);
  \node[](a1)at(2.7cm,0){\bf 2};
    \draw[-,line width=0.15em,color=blue,fill=ballblue] (0.835,0) arc (0:180:0.3);
\draw[-,line width=0.15em,color=red,fill=red!30] (0.23,0) arc (180:360:0.3);
  \node[](a2)at(0.523cm,0){\bf 1};
\node[bluenode](a3)[xshift=2.15cm,yshift=-1.7cm]{\bf 3}; 
 \draw[black,<->,>=stealth](0.85,0)--(2.35,0);
\draw[black,->,>=stealth](0.6,-0.27)--(1.5,-1.12);
\draw[black,->,>=stealth](2.7,-0.28)--(1.9,-1.08);
\draw[->,>=stealth,black](a2)edge[in=-20,out=60,loop above]node[below right]{}();
\draw[->,>=stealth,black](a1)edge[in=-20,out=60,loop above]node[below right]{}();
\end{tikzpicture}
\end{figure}
 \end{example}

\begin{example}\label{ex215}\rm  Theorem \ref{thm:BotNck} imposes graph-theoretic restrictions on systems that can have robust solutions.
 In the  bi-directional graph in Figure~\ref{fg:4-troph}(a), each blue node species on each level is connected {\it only} to the species in the adjacent level, though some edges connect red nodes with other red nodes.
Let $N_1, N_2, N_3,N_4$
denote the number of species in each of these ``trophic'' levels, listing from the bottom to the top. {\sana This type of trophic graph cannot have robust solutions if the total number of species in the odd-numbered levels, $N_{\mbox{odd}} = N_1+N_3$ is greater than the total number of species in the even-numbered levels. For example, in Figure~\ref{fg:4-troph}(a),
$N_{\mbox{odd}}= 8 > N_{\mbox{even}}=5$, so there is a 3-backward bottleneck $\{{\rm{red\ nodes}\} \rightarrow\{\rm blue\ nodes}\}.$
By the Bottleneck Theorem~\ref{thm:BotNck}, at least 3 edges must be added to allow robust solutions.}
\end{example}

In Appendix B, by using a numerical approach we identify the set of kernel nodes for this $N=M=13$ system. One finds that the set of kernel nodes consists of the eight blue nodes, and therefore they constitute the bottle of the backward bottleneck (and also the forward bottleneck $\{{\rm blue\ nodes}\}\rightarrow \{{\rm red\ nodes}\}$, since the graph is symmetric). That is, for almost every $13\times 13$ Jacobian matrix of the structured system shown in Fig.~\ref{fg:4-troph}(a), all vectors in the three-dimensional nullspace have zero entries for the coordinates corresponding to the red nodes.

\begin{figure}
\begin{center}
\begin{subfigure}[b]{.4\textwidth}
\centering
\begin{tikzpicture}
    [ultra thick]
    \node[rdnode](a11)[xshift=-.5cm,yshift=1cm]{};
    \node[rdnode](a12)[xshift=+.5cm,yshift=1cm]{};
   \node[bdnode](a21)[xshift=-1.5cm,yshift=0cm]{};
    \node[bdnode](a22)[xshift=-.5cm,yshift=0cm]{};
    \node[bdnode](a23)[xshift=+.5cm,yshift=0cm]{};
    \node[bdnode](a24)[xshift=+1.5cm,yshift=0cm]{};
    \node[rdnode](a31)[xshift=-1.4cm,yshift=-1cm]{};
    \node[rdnode](a32)[xshift=0cm,yshift=-1cm]{};
    \node[rdnode](a33)[xshift=1.4cm,yshift=-1cm]{};
   \node[bdnode](a41)[xshift=-1.8cm,yshift=-2cm]{};
    \node[bdnode](a42)[xshift=-.5cm,yshift=-2cm]{};
    \node[bdnode](a43)[xshift=+.5cm,yshift=-2cm]{};
    \node[bdnode](a44)[xshift=+1.8cm,yshift=-2cm]{};
\draw[<->,>=stealth,black] (a41)-- (a32);
\draw[<->,>=stealth,black] (a44)-- (a33);
\draw[<->,>=stealth,black] (a44)-- (a31);
\draw[<->,>=stealth,black] (a43)-- (a32);
\draw[<->,>=stealth,black] (a43)-- (a33);
\draw[<->,>=stealth,black] (a42)-- (a31);
\draw[<->,>=stealth,black] (a42)-- (a33);
\draw[<->,>=stealth,black] (a41)-- (a32);
\draw[<->,>=stealth,black] (a41)-- (a33);
\draw[<->,>=stealth,black] (a31)-- (a21);
\draw[<->,>=stealth,black] (a31)-- (a22);
\draw[<->,>=stealth,black] (a31)-- (a24);
\draw[<->,>=stealth,black] (a32)-- (a21);
\draw[<->,>=stealth,black] (a32)-- (a22);
\draw[<->,>=stealth,black] (a33)-- (a23);
\draw[<->,>=stealth,black] (a33)-- (a21);
\draw[<->,>=stealth,black] (a32)-- (a24);
\draw[<->,>=stealth,black] (a21)-- (a11);
\draw[<->,>=stealth,black] (a21)-- (a12);
\draw[<->,>=stealth,black] (a22)-- (a11);
\draw[<->,>=stealth,black] (a22)-- (a12);
\draw[<->,>=stealth,black] (a23)-- (a11);
\draw[<->,>=stealth,black] (a23)-- (a12);
\draw[<->,>=stealth,black] (a24)-- (a11);
\draw[<->,>=stealth,black] (a24)-- (a12);
\draw[<->,>=stealth,black](a31)to[bend right,out=120,in=105](a11);
\draw[<->,>=stealth,black] (a33)to[bend right,out=-100,in=-100](a12);
\node at ([xshift=.0cm,yshift=.5cm]a11.north east){\Huge{\xmark}};
\node at ([xshift=.5cm,yshift=.5cm]a11.north east){\Huge{\xmark}};
\node at ([xshift=1.cm,yshift=.5cm]a11.north east){\Huge{\xmark}};
\end{tikzpicture}
\caption{}
\label{fg:4-troph d}
\end{subfigure}
\begin{subfigure}[b]{.45\textwidth}
     \centering
     \resizebox{4cm}{4cm}{%
     \begin{tikzpicture}[ultra thick]
       \tikzset{every node}=[font=\bf]
\node[bluenode](s11)[xshift=10.5cm,yshift=0.5cm]{1};
\node[bluenode](s12)[xshift=12cm,yshift=0.5cm]{2};
\node[bluenode](s13)[xshift=14cm,yshift=0.5cm]{3};
\node[bluenode](s14)[xshift=15.5cm,yshift=0.5cm]{4};
  \draw[-,line width=0.15em,color=blue,fill=ballblue] (12.025,-1.2) arc (0:180:0.33);
\draw[-,line width=0.15em,color=red,fill=red!30] (11.37,-1.2) arc (180:360:0.33);
  \node[](s22)at(11.7cm,-1.2cm){6};
    \draw[-,line width=0.15em,color=blue,fill=ballblue] (9.34,-1.2) arc (0:180:0.33);
\draw[-,line width=0.15em,color=red,fill=red!30] (8.687,-1.2) arc (180:360:0.33);
  \node[](s21)at(9.023cm,-1.2cm){5};
\node[rednode](s31)[xshift=13.1cm,yshift=-3.3cm]{7};
 \draw[<-,>=stealth,black](8.3cm,1.5cm)-- (s11);
 \draw[<-,>=stealth,black](9.6cm,1.5cm)-- (s12);
 \draw[<-,>=stealth,black](11.3cm,1.5cm)-- (s13);
\draw[<-,>=stealth,black](12.5cm,1.5cm)-- (s14);
  \draw[->,>=stealth,black](s21)-- (s12);
    \draw[->,>=stealth,black](s21)-- (s11);
    \draw[->,>=stealth,black](s22)-- (s13);
    \draw[->,>=stealth,black](s22)-- (s14);  
     \draw[<-,>=stealth,black](s21)-- (s31);
    \draw[<-,>=stealth,black](s22)-- (s31);
     \draw[<->,>=stealth,black](9cm,-2.6)-- (s31);
     \node at ([xshift=1cm,yshift=1.2cm]s11.north east){\Huge{\xmark}};
\node at ([xshift=2cm,yshift=1.2cm]s11.north east){\Huge{\xmark}};
\node at ([xshift=3cm,yshift=1.2cm]s11.north east){\Huge{\xmark}};
    \end{tikzpicture}}
    \caption{}
    \label{fg:4-troph e}
\end{subfigure}
\end{center}
\caption{
{\bf Examples of robust and fragile graphs.} 
A red X-mark indicates that every $F$ that respects the graph is fragile; such an $F$ has no robust steady states.
The number of  X-marks is the minimum number of edges that must be added before the graph can be robust. 
{\bf (a)} The graph discussed in Example \ref{ex215} is fragile. The three X-marks mean there is a backward 3-bottleneck.
{\bf  (b)} 
Illustration of how bottlenecks combine to form the minimax bottleneck. A portion of a larger network is shown.
For a backward bottleneck, edges leaving the nodes of the bottle to nodes in the rest of the graph can occur without changing the bottleneck, provided they do not point to bottleneck nodes.
Here, the sets $\{1,2\}$, $\{3,4\}$, $\{1,2,3,4\}$, and $\{5,6\}$ each form the bottle of a distinct backward bottleneck. Combining them, $B = \{1, 2, 3, 4, 5, 6\}$ is the
bottle of the  \minimax 3-bottleneck that contains the smaller bottlenecks. The neck of the minimax backward bottleneck is $B^{\fm} = \{5, 6, 7\}$.}
\label{fg:4-troph}
 \end{figure}
\begin{example}\rm  \label{ex26nodes}
 Figure~\ref{fg:26}(a) analyzes Figure~\ref{fg:12}.
Our results show that the graph is fragile. There is a 6-species bottleneck. In order to eliminate the bottleneck by adding edges, six edges must be carefully added.  A substantial research effort has aimed at discovering what stabilizes ecological networks. Some have described ways that appear to promote the stability of an ecosystem \cite{gross_generalized_2009,maynard_phenotypic_2019}, but there has been no gold standard for assessing robustness, which is a requirement for stability.
Gross et al.~\cite{gross_generalized_2009} suggested two universal rules:
Food-web stability is enhanced when (i) ``species at a high trophic level feed on multiple prey species'',
and (ii) ``species at an intermediate trophic level are fed upon by multiple predator species.''
These rules suggest adding edges might enhance the stability of the ecological network, but it is not the full story. As an example, 
Figure~\ref{fg:26}(b) is made by adding 29 bi-directional (green) edges to Figure~\ref{fg:26}(a),
adding edges throughout the network, but the resulting network, Figure~\ref{fg:26}(b) is not yet robust. Our theory suggests where the edges must be added to make a system robust.
Write $B$ for the set of blue nodes and
$B^\fm$ the red nodes. Since $\#(B)=12$ and $\#(B^\fm)=6$, $B$ is a 6-bottleneck so the graph cannot be robust. At least six backward edges ending in $B$ and not starting in $B^\fm$ must be added before the graph becomes robust.
\end{example}
\begin{figure}
\begin{center}
\begin{subfigure}[c]{\columnwidth}
\centering
\resizebox{11cm}{4.5cm}{%
 \begin{tikzpicture}
    [ultra thick]
    \tikzset{every node}=[font=\BF]
    \node[2uninode]at(-5.7cm,2cm)(a11){$1$};
    \node[2uninode]at(-3.5cm,2cm)(a12){${2}$};
    \node[2uninode]at(0cm,2cm)(a13){${3}$};
    \node[2secnode]at(2.2cm,2cm)(a14){${4}$};
    \node[2uninode]at(5cm,2cm)(a15){${5}$};
    \node[2tednode]at(+7cm,2cm)(a16){${6}$};
    \node[2tednode]at(+9cm,2cm)(a17){${7}$};
    \node[2uninode]at(-6.5cm,.2cm)(a21){${8}$};
    \node[2uninode]at(-4.8cm,.2cm)(a22){${9}$};
    \node[tednode]at(-3cm,.2cm)(a23){${10}$};
    \node[uninode]at(-0.7cm,.2cm)(a24){${11}$};
    \node[tednode]at(1.5cm,.2cm)(a25){${12}$};
    \node[uninode]at(3.5cm,.2cm)(a26){${13}$};
    \node[tednode]at(5.5cm,.2cm)(a27){${14}$};
    \node[secnode]at(+8cm,.2cm)(a28){${15}$};
     \node[secnode]at(-5cm,-1.5cm)(a31){${16}$};
    \node[secnode]at(-1.5cm,-1.5cm)(a32){${17}$};
    \node[secnode]at(3cm,-1.5cm)(a33){${18}$};
    \node[secnode]at(6.5cm,-1.5cm)(a34){${19}$};
    \node[tednode]at(-6cm,-3.2cm)(a41){${20}$};
    \node[tednode]at(-4cm,-3.2cm)(a42){${21}$};
    \node[tednode]at(-1cm,-3.2cm)(a43){${22}$};
    \node[tednode]at(1cm,-3.2cm)(a44){${23}$};
    \node[tednode]at(+4cm,-3.2cm)(a45){${24}$};
    \node[tednode]at(6cm,-3.2cm)(a46){${25}$};
    \node[tednode]at(8cm,-3.2cm)(a47){${26}$};
   \draw[<->,>=stealth,black] (a11)-- (a21);
    \draw[<->,>=stealth,black] (a12)-- (a22);
    \draw[<->,>=stealth,black] (a13)-- (a24);
    \draw[<->,>=stealth,black] (a14)-- (a25);
    \draw[<->,>=stealth,black] (a14)-- (a27);
    \draw[<->,>=stealth,black] (a14)-- (a33);
    \draw[<->,>=stealth,black] (a14)-- (a33);
    \draw[<->,>=stealth,black] (a15)-- (a26);
    \draw[<->,>=stealth,black] (a16)-- (a28);
    \draw[<->,>=stealth,black] (a17)-- (a28);
    \draw[<->,>=stealth,black] (a21)-- (a31);
    \draw[<->,>=stealth,black] (a21)-- (a32);
     \draw[<->,>=stealth,black] (a22)-- (a32);
    \draw[<->,>=stealth,black] (a23)-- (a31);
     \draw[<->,>=stealth,black] (a23)-- (a32);
     \draw[<->,>=stealth,black] (a23)-- (a33);
     \draw[<->,>=stealth,black] (a23)-- (a34); 
    \draw[<->,>=stealth,black] (a24)-- (a32);
    \draw[<->,>=stealth,black] (a25)-- (a32);
    \draw[<->,>=stealth,black] (a26)-- (a33);
    \draw[<->,>=stealth,black] (a27)-- (a32);
 \draw[<->,>=stealth,black] (a27)-- (a34);
  \draw[<->,>=stealth,black] (a28)-- (a34);
  \draw[<->,>=stealth,black] (a34)-- (a47);
  \draw[<->,>=stealth,black] (a34)-- (a46);
  \draw[<->,>=stealth,black] (a34)-- (a45);
  \draw[<->,>=stealth,black] (a34)-- (a44);
  \draw[<->,>=stealth,black] (a33)-- (a44);
  \draw[<->,>=stealth,black] (a31)-- (a44);
  \draw[<->,>=stealth,black] (a31)-- (a42);
  \draw[<->,>=stealth,black] (a32)-- (a42);
  \draw[<->,>=stealth,black] (a32)-- (a41);
  \draw[<->,>=stealth,black] (a32)-- (a43);
  \draw[<->,>=stealth,black] (a32)-- (a44);
\node at ([xshift=5.0cm,yshift=.5cm]a11.north east){\Huge{\xmark}};
\node at ([xshift=5.5cm,yshift=.5cm]a11.north east){\Huge{\xmark}};
\node at ([xshift=6.0cm,yshift=.5cm]a11.north east){\Huge{\xmark}};
\node at ([xshift=6.5cm,yshift=.5cm]a11.north east){\Huge{\xmark}};
\node at ([xshift=7.0cm,yshift=.5cm]a11.north east){\Huge{\xmark}};
\node at ([xshift=7.5cm,yshift=.5cm]a11.north east){\Huge{\xmark}};
\end{tikzpicture}
}
\caption{}
\label{fg:26 a}
\end{subfigure}
\begin{subfigure}[c]{\columnwidth}
\centering
\resizebox{11cm}{4.5cm}{%
 \begin{tikzpicture}
    [ultra thick]
     \tikzset{every node}=[font=\BF]
      \node[2uninode]at(-5.7cm,2cm)(a11){${1}$};
    \node[2uninode]at(-3.5cm,2cm)(a12){${2}$};
    \node[2uninode]at(0cm,2cm)(a13){${3}$};
    \node[2secnode]at(2.2cm,2cm)(a14){${4}$};
    \node[2uninode]at(5cm,2cm)(a15){${5}$};
    \node[2tednode]at(+7cm,2cm)(a16){${6}$};
    \node[2tednode]at(+9cm,2cm)(a17){${7}$};
    \node[2uninode]at(-6.5cm,.2cm)(a21){${8}$};
    \node[2uninode]at(-4.8cm,.2cm)(a22){${9}$};
    \node[tednode]at(-3cm,.2cm)(a23){${10}$};
    \node[uninode]at(-0.7cm,.2cm)(a24){${11}$};
    \node[tednode]at(1.5cm,.2cm)(a25){${12}$};
    \node[uninode]at(3.5cm,.2cm)(a26){${13}$};
    \node[tednode]at(5.5cm,.2cm)(a27){${14}$};
    \node[secnode]at(+8cm,.2cm)(a28){${15}$};
     \node[secnode]at(-5cm,-1.5cm)(a31){${16}$};
    \node[secnode]at(-1.5cm,-1.5cm)(a32){${17}$};
    \node[secnode]at(3cm,-1.5cm)(a33){${18}$};
    \node[secnode]at(6.5cm,-1.5cm)(a34){${19}$};
    \node[tednode]at(-6cm,-3.2cm)(a41){${20}$};
    \node[tednode]at(-4cm,-3.2cm)(a42){${21}$};
    \node[tednode]at(-1cm,-3.2cm)(a43){${22}$};
    \node[tednode]at(1cm,-3.2cm)(a44){${23}$};
    \node[tednode]at(+4cm,-3.2cm)(a45){${24}$};
    \node[tednode]at(6cm,-3.2cm)(a46){${25}$};
    \node[tednode]at(8cm,-3.2cm)(a47){${26}$};
    \draw[<->,>=stealth,black] (a11)-- (a21);
    \draw[<->,>=stealth,black] (a12)-- (a22);
    \draw[<->,>=stealth,black] (a13)-- (a24);
    \draw[<->,>=stealth,black] (a14)-- (a25);
    \draw[<->,>=stealth,black] (a14)-- (a27);
    \draw[<->,>=stealth,black] (a14)-- (a33);
    \draw[<->,>=stealth,black] (a14)-- (a33);
    \draw[<->,>=stealth,black] (a15)-- (a26);
    \draw[<->,>=stealth,black] (a16)-- (a28);
    \draw[<->,>=stealth,black] (a17)-- (a28);
    \draw[<->,>=stealth,black] (a21)-- (a31);
    \draw[<->,>=stealth,black] (a21)-- (a32);
     \draw[<->,>=stealth,black] (a22)-- (a32);
    \draw[<->,>=stealth,black] (a23)-- (a31);
     \draw[<->,>=stealth,black] (a23)-- (a32);
     \draw[<->,>=stealth,black] (a23)-- (a33);
     \draw[<->,>=stealth,black] (a23)-- (a34); 
    \draw[<->,>=stealth,black] (a24)-- (a32);
    \draw[<->,>=stealth,black] (a25)-- (a32);
    \draw[<->,>=stealth,black] (a26)-- (a33);
    \draw[<->,>=stealth,black] (a27)-- (a32);
 \draw[<->,>=stealth,black] (a27)-- (a34);
  \draw[<->,>=stealth,black] (a28)-- (a34);
  \draw[<->,>=stealth,black] (a34)-- (a47);
  \draw[<->,>=stealth,black] (a34)-- (a46);
  \draw[<->,>=stealth,black] (a34)-- (a45);
  \draw[<->,>=stealth,black] (a34)-- (a44);
  \draw[<->,>=stealth,black] (a33)-- (a44);
  \draw[<->,>=stealth,black] (a31)-- (a44);
  \draw[<->,>=stealth,black] (a31)-- (a42);
  \draw[<->,>=stealth,black] (a32)-- (a42);
  \draw[<->,>=stealth,black] (a32)-- (a41);
  \draw[<->,>=stealth,black] (a32)-- (a43);
  \draw[<->,>=stealth,black] (a32)-- (a44);
  \draw[<->,>=stealth,bgreen](a17)to [bend right,out=25,in=150](a34);
  \draw[<->,>=stealth,bgreen](a28)to [bend right,out=370,in=150](a47);
  \draw[<->,>=stealth,bgreen](a14)-- (a24);
   \draw[<->,>=stealth,bgreen](a31)-- (a32);
   \draw[<->,>=stealth,bgreen](a11)-- (a12);
  \draw[<->,>=stealth,bgreen](a12)to [bend right,out=10,in=170](a32);
  \draw[<->,>=stealth,bgreen](a14)to [bend right,out=10,in=170](a44);
  \draw[<->,>=stealth,bgreen](a16)to [bend right,out=-30,in=150](a33);
   \draw[<->,>=stealth,bgreen] (a14)-- (a32);
  \draw[<->,>=stealth,bgreen] (a14)-- (a23);
  \draw[<->,>=stealth,bgreen](a15)to [bend right,out=0,in=170](a33);
 \draw[<->,>=stealth,bgreen] (a13)-- (a22);
 \draw[<->,>=stealth,bgreen] (a13)-- (a14);
\draw[<->,>=stealth,bgreen] (a14)-- (a15);
\draw[<->,>=stealth,bgreen] (a28)-- (a44);
 \draw[<->,>=stealth,bgreen] (a11)-- (a22);
 \draw[<->,>=stealth,bgreen] (a46)-- (a33);
 \draw[<->,>=stealth,bgreen] (a45)-- (a33);
 \draw[<->,>=stealth,bgreen] (a43)-- (a31);
 \draw[<->,>=stealth,bgreen] (a27)-- (a28);
\draw[<->,>=stealth,bgreen] (a33)-- (a42);
\draw[<->,>=stealth,bgreen] (a25)-- (a34);
\draw[<->,>=stealth,bgreen] (a26)-- (a34);
\draw[<->,>=stealth,bgreen] (a41)-- (a31);
\draw[<->,>=stealth,bgreen] (a32)-- (a33);
\draw[<->,>=stealth,bgreen] (a12)-- (a13);
\draw[<->,>=stealth,bgreen] (a21)-- (a22);
\draw[<->,>=stealth,bgreen] (a12)-- (a24);
\draw[<->,>=stealth,bgreen] (a22)-- (a31);
\node at ([xshift=5.0cm,yshift=.5cm]a11.north east){\Huge{\xmark}};
\node at ([xshift=5.5cm,yshift=.5cm]a11.north east){\Huge{\xmark}};
\node at ([xshift=6.0cm,yshift=.5cm]a11.north east){\Huge{\xmark}};
\node at ([xshift=6.5cm,yshift=.5cm]a11.north east){\Huge{\xmark}};
\node at ([xshift=7.0cm,yshift=.5cm]a11.north east){\Huge{\xmark}};
\node at ([xshift=7.5cm,yshift=.5cm]a11.north east){\Huge{\xmark}};
\end{tikzpicture}}
\caption{}
\label{fg:26 b}
\end{subfigure}
\caption{{\bf { Existence of a backward} bottleneck makes a network fragile.} The six $X$'s mean six species must be eliminated before it becomes robust -- unless carefully chosen edges are added. 
 (a) 
 Let $B$ denote the 12 nodes shaded blue. Then the backward set $B^\fm$ consists of 6 nodes shaded red. Hence $B$ is a 6-species bottleneck. At least six edges that increase the set $B^\fm$ 
must be added to make the graph robust. 
It is easy to cover 20 nodes with disjoint cycles, leaving 6 in $B$ uncovered, so the structure matrix has rank 20.
(b) 
29 bi-directional edges have been added, shown in green, adding at least one edge to every node, and edges between every pair of layers, and edges within each layer except the lowest. These make no difference to the robustness. $B$ is still a 6-species bottleneck. 
}\label{fg:26}
\end{center}
\end{figure}
\section{Discussion}\label{s:discussion}

{The Competitive Exclusion Principle (CEP), long discussed in ecological literature,
states that two predators that depend solely on the same prey species cannot coexist, as portrayed in Figure~\ref{fg:Fragile}. More precisely, the CEP says the two predators whose population density depends purely on a prey species population density cannot coexist unless they benefit precisely equally from the prey, a scenario that is virtually impossible in natural circumstances.} {\jim In terms of the graph, there are two ``predator'' nodes that depend upon one ``prey'' node. That means there are two nodes for which each has only one incoming edge and two edges come from the same node, a third node. This simple configuration is a bottleneck. Since there is a backward  bottleneck,  by Theorem~\ref{thm:BotNck}(III), the graph is fragile (\ie, not robust). 
}

We show in this article how the basic lesson of the Competitive Exclusion Principle of ecology can be extended to a general concept that applies to all systems $F(x)=c$ of $M$ equations in $N$ unknowns, provided that they respect an underlying structure, a structure that can be encapsulated in a vector space of functions $F$.

 If there is a backward bottleneck, the bottle can be viewed as being a set of predator species and the neck as prey species. 
The bottle consists of $n$ species that are only influenced by $n-k$ species where $k>0$; 
 It does not matter whether the influence is positive (such as a prey or food source) or is negative (such as a predator). The existence of a bottleneck can be detected by looking at the rank of the whole system. Ihere is a $k$-backward bottleneck somewhere in the system if and only if  the rank $r$ of the system satisfies $r \leq M-k$. The same is true for forward bottlenecks with $M$ replaced by $N$.

 In Theorem \ref{thm:BotNck}, we proved the equivalence of two ostensibly different views of structured systems, the rank of the structure matrix, and the existence of bottlenecks. 
Depending on what is known about the network model, or class of models, one or another of the views may be most informative. 
These results have immediate implications to possible graph structures, for example precluding robust solutions in strictly trophic food webs without interactions between competing species unless strict constraints on the number of species in each guild are satisfied. 
These constraints are direct generalizations of the CEP. Interestingly, they may have extensive implications for systems in general, outside ecology. 

When $M=N$, the existence of a forward bottleneck implies the existence of a backward bottleneck, and vice versa. When the generic rank of a structured matrix is strictly less than $M=N$, both bottlenecks exist and the system is fragile.

 When $N \neq M$, the existence of forward and backward bottlenecks are more independent. Let the maxrank of the $M$ by $N$ structure matrix be $r$:
 \begin{itemize}
\item  if a collection of  $M$ equations depends upon $M-r$ fewer variables, then there exists a unique minimax backward bottleneck with rank deficiency of $M-r$. Hence, for  almost every $F$ that respects the graph $G$ there are no robust solutions.

\item   if a collection of $r$ equations are overly determined by $N$ variables where $N-r>0$, then there exists a unique minimax forward bottleneck with rank deficiency of $N-r$.  Then, for \textit{almost} every $F$ that respects the graph $G$ and for almost every $x$, ${\cSx}$ is a manifold of dimension $N-r$.
\end{itemize}
In particular, when $N\neq M$, existence of a forward bottleneck does not rule out the existence of robust solutions.

From our perspective, the Competitive Exclusion Principle (CEP) is not a biological principle whose {\it conclusion} can be tested. It is a theorem whose {\it hypotheses} can be tested in various biological settings. It is a theorem that we have generalized here, in the form of Theorems~\ref{thm:BotNck} and \ref{thm:CR}.  These ideas about equations and networks can be extended far beyond the ecological realm. 
 
The Bottleneck Theorem says the correct generalization of the Competitive Exclusion Principle is not simply about predators and prey. It says that if a collection of species depends purely on $k$ fewer species, where $k>0$, there can not be a robust steady state. That is a statement about the equations of nature. This fact generalizes to a statement about the nature of structured equations: If a collection of $n$ equations depends upon $n-k $ variables where $k >0$, then there exists a unique minimax backward $k$-bottleneck. Hence, for  almost every $F$ that respects the graph $G$, there are no robust solutions of $F(x)=c$ for any c, and for almost every $c$, and
almost every $F$ that respects the graph $G$, there are no solutions, not even fragile solutions.
{\bf The Bottleneck Theorem is an extended version of the Competitive Exclusion Principle, applicable to scientific and engineering areas well beyond ecology.}

\section*{Acknowledgements}
 Thanks to Dima Dolgopyat,  Shuddho Das, and Roberto De Leo for their helpful comments.
 SJ was partially supported by CIHR and NSERC through grants to James Watmough (RGPIN-2017-05760). TS was partially supported by the NSF (grant DMS1723175). We thank the reviewers for several helpful comments that led to improvements in the manuscript.
\section*{Author contributions statement}
All the authors contributed equally to this project.

\newpage
\section*{Appendix A}\label{s:appendix}
We begin with some lemmas which support the definition of kernel nodes. The proof of Theorem \ref{thm:BotNck} follows.

\begin{lemma} \label{newlem1}
Let $w, v_1, \ldots, v_k \in R^n$. Then $w$ is in Span$\{v_1,\ldots, v_k\}$ if and only if there is a subset $\{v_{i_1},\ldots,v_{i_r}\}$ such that rank $\{v_{i_1},\ldots,v_{i_r}\}$ = r = rank $\{w, v_{i_1},\ldots,v_{i_r}\}$
\end{lemma}

\begin{lemma}\label{newlem2}
If $S_1,\ldots,S_k$ are Lebesgue-measurable subsets of $R^n$ of either full measure or measure zero, then the same is true for all finite unions and intersections of the $S_i$ and their complements.
\end{lemma}

Consider any vector space of $m\times n$ matrices, such as the space $\cL(G)$ for some directed graph $G$.
Assume there is a prescribed set of matrix minors.
The next lemma states that the set of matrices 
on which those matrix minors all vanish either
has measure zero or has
full measure -- in the Lebesgue measure on the vector space.

\begin{lemma}\label{newlem3}
Let $A^i$ be an $m\times n$ matrix for $i = 1,\ldots,k$.
Consider the $k$-dimensional set of parametrized matrices $C = c_1A^1+\ldots + c_kA^k$ (where the $c_j$ are the parameters). Let
$M_1(c_1, \ldots,c_k), \ldots,\  M_m(c_1, \ldots, c_k)$ enumerate all minor determinants of $C$. Then for any subset $S$ of the integers $\{1, \ldots, m\}$, the set of $(c_1,\ldots,c_k)$ for which $M_i(c_1,\ldots,c_k) = 0$ if and only if $i\in S$ has either full measure or measure zero in $R^k$.
\end{lemma}
\begin{proof}
    
 All the $M_i$ are polynomials in the $c_i$ and are either identically zero (vanish on a full measure set) or not identically zero (vanish on a lower dimensional, thus measure zero set). Specifying some of the $M_i$ to be zero and the rest to be nonzero results in a full measure or zero measure subset of $(c_1,\ldots,c_k)$, according to Lemma \ref{newlem2}.
 \end{proof}

\begin{lemma}\label{newlem4}
Let $H$ be a dimension $n-1$ subspace of $R^n$ and $C = c_1A^1+\ldots + c_kA^k$. Then the subset of $(c_1,\ldots,c_k)\in R^k$ for which ker$(C) \subset H$ is either full measure or measure zero.
\end{lemma}
\begin{proof}
 Without loss of generality, we may assume $H = \{x_1 = 0\}$. If the first column $C_1$ of $C$ is identically zero, then the subset has full measure. If not, then ker$(C) \not\subset H$ for some $(c_1,\ldots,c_k)$ if and only if the column $C_1$ lies in the span of the rest of the columns of $C$. By Lemma \ref{newlem1}, this holds exactly for $(c_1,\ldots,c_k)$ in a subset defined by minor determinants, which has either full measure or zero measure by Lemma \ref{newlem3}.
\end{proof}

\begin{lemma} \label{lm20}
For each integer $1 \leq j \leq N$, either $S_j$ or its complement $S_j^c$ has measure zero in $\cL(G)$.
\end{lemma}
\begin{proof}
 Note that $S_j$ is the set of matrices with structure matrix $S(G)$ for which the kernel is contained in $H = \{x_j = 0\}$. Lemma \ref{newlem4} says that either this set or its complement has measure zero.
\end{proof}
Lemma~\ref{lm20} can be illustrated explicitly for $2\times 2$ structure matrices $A$.
Let 
$$A = \left[\begin{array}{cc}
  a_{11}  & a_{12} \\
a_{21}     & a_{22}
\end{array}\right],$$ where some of the entries may be required to be 0 in the structure matrix. 
Consider the set of matrices
$$S_1 = \{A|\  a_{11}a_{22}-a_{12}a_{21} \neq 0\} \cup \{A|\  a_{12}=0  \text{ and  } a_{11}\neq 0\} \cup 
\{A|\  a_{22}=0 \text{ and } a_{21}\neq 0 \},
$$
and its complement
$$S_1^c =\{A|\ a_{11}a_{22}-a_{12}a_{21}= 0\} \cap \{A|\ a_{12}\neq 0  \text{ or } a_{11}= 0\} \cap 
\{A|\ a_{22}\neq 0 \text{ or } a_{21} = 0 \}.
$$
If the structure matrix allows all entries of $A$ to be nonzero, then $(S_1^c)$  is of measure zero, and so $S_1$ has full measure. 

As another example, if the structure matrix requires $a_{11}=a_{21}=0$, the roles are reversed. In that case $S_1$ is of measure zero 
and so $S_1^c$ has full measure. 

\bigskip

The proof of Theorem~\ref{thm:BotNck} follows.
\begin{proof}[Proof of Theorem~\ref{thm:BotNck} Part (I)]
{\sana Assume the pair of sets $(B, B^\fm)$ is a backward $K$-bottleneck such that $K>0$, $B$ and $B^\fm$ have $b$ and $b-K$ nodes, respectively, where $b\geq 2$.} That means the $b$ rows of the matrix $S$ corresponding to $B$ are zero except for $b-K$ columns, or in other words, $S$ has form
\begin{equation} \label{Smat}
S = 
\overset{\ \Huge{N-b+K}\ \Large{b-K}}{
\left[
\begin{array}{c|c}
\ \ \ 0 \ \ & \star\\ \hline  \ \ \star & \star\\
\end{array}
\right]
}
\begin{array}{l}
 b \\ M-b
\end{array}\ \ \ \ \ \ \ \ \ 
T = 
\overset{\ \Huge{b} \   \Large{ M-b}}{
\left[
\begin{array}{c|c}
0 & \star\\ \hline  \star & \star\\
\end{array}
\right]
}
\begin{array}{l}
{N}-b+K \\ b-K
\end{array}
\end{equation}
where we have renumbered  the nodes so that $B$ consists of the first $b$ rows of $S$  and where we have situated 
the $b-K$ nonzero columns at the right of the matrix for simplicity (the nonzero columns could be anywhere in $S$).
The transpose $T$ is also shown for convenience.

It follows that the  $N \times b$ submatrix $T_B$ of $T$ consisting of the $b$ columns represented by the nodes of $B$ has at most $b-\M$ nonzero rows, the rows corresponding to nodes of $B^\fm$. Since $T_B:R^b\to R^N$ has at most $b-K$ nonzero rows, the kernel of $T_B$ is at least $\M$-dimensional, which implies the same about the kernel of $T$, and so $r = $ rank$(T) \leq M-K<M$, as required. 

To prove the converse direction, we will assume $r<M$, and show that 
$S$ must have form (\ref{Smat}) with $K=M-r$, with the proviso, as above, that the rightmost $b-K$ nonzero columns could occur anywhere in $S$. 
If that can be shown, then there is a backward $K$-bottleneck,
with $K = M-r > 0$.

Let $B = B_{\rm kernel}$ be the set of kernel nodes of $T$. According to Lemma~\ref{lm20}, for any node $j$ in $B$, the set $T_j$ of $N\times M$ matrices respecting $T$ for which $x_j=0$ for every $x\in\ $ker $(T)$ has measure zero. Let $b=|B|$.

First, note that if $B$ is the empty set, then the intersection of the {$T_j$} is a full measure set, i.e. for almost every $A$ respecting the structure {$T$}, $\ker(A) = \{0\}$, a contradiction to 
$r<M$.
Thus $b>0$. 
Note that if $j$ is in the complement $B^c$, then $x_j = 0$ for every vector $x$ in $\ker(A)$, for every $A$ not in the measure zero set {$Q = \bigcup_{j\in B^c} T_j^c$.}  We use this fact below.

Consider matrices $A$ in $\bm{\cL(G)}$ that respect the structure matrix $T$, and are not in the set {$Q$}. Renumber the nodes such that $B = \left\{1, \ldots, b\right\}.$  Let $A_B$ denote the submatrix of first $b$ columns of $A$.
 Since $A$ is not in {$Q$}, $x_{b+1}=\ldots =x_{M}=0$ for $x$ in ker$(A)$, so we can assume ker$(A)=  [U, 0, \ldots, 0]$, where $U$ is a dimension $K$  subspace of $\R^b$ (since the rank of $A$ is $r=M-K$), and the space $U^{\perp}$ is $b-K$ dimensional. 
For almost every $A$ respecting {$T$, the rows $r_1, \ldots, r_N$} of the submatrix $A_B$ must all be in $U^{\perp}$, and satisfy the following two properties, proved below: (1) the nonzero rows of $A_B$ are linearly independent, and (2) no set of $p$ columns of $A_B$ contains all of the nonzero entries of $p$ or more nonzero rows. 
Property (1)  forces all but $b-K$ of the rows of $A_B$ to be zero rows in the structure $T$. Thus the first $b$ columns of $A$ have at most $b-K$ nonzero rows, which verifies the form of $T$ in (\ref{Smat}).

Finally, we verify (1) and (2). If $b=1$, (1) is true because if any entry in the (single) column of $A_B$ is nonzero, the first component of vectors in $\ker(A)$ is zero for almost every $A$, a contradiction to the definition of $B$. For $b>1$, we induct on $b$. Let $0 = \sum c_ir_i$ be a dependency of rows of $A_B$, where all $c_i\neq 0$. If the union of the coordinates appearing in the rows of the $r_i$ does not include all $b$ coordinates, use the induction hypothesis. If they include all $b$ coordinates, then so does the structure matrix $T$, and no such dependency can exist for almost every $A$. To verify (2), suppose there are columns $c_{i_1},\ldots,c_{i_p}$ and $p$ such rows. By (1), the rows are linearly independent, and their entries are restricted to $p$ columns. Therefore, any vector $x=(u,0,\ldots,0)$ in ker$(A)$ must be zero in entries $u_{i_1},\ldots,u_{i_p}$. This contradicts the fact that $i_1,\ldots,i_p$ are kernel nodes. 
\end{proof}

\begin{proof}[Proof of Theorem~\ref{thm:BotNck}, Part (II)]  Note that $B_{\rm kernel}$ is the $K^*$-bottleneck used\\ above in the first part of the proof, where $K^*$ is such that {rank$(S)$=rank$(T) = M-K^*$}. No bottleneck exists with larger $K$, by Part (I). Also, $B$ is minimal because no node can be deleted without the bottleneck becoming a $K$-bottleneck for $K<K^*$, due to property (2) above.  Therefore $B_{\rm kernel}$ of $T$ is the bottle of the \minimax backward bottleneck of $S$.
\end{proof}

 \begin{proof}[Proof of Theorem~\ref{thm:BotNck}, Part (III)]
Part (I) of this theorem implies that there exists a backward $K$-bottleneck if and only if $M-r \geq K>0$. Hence, there exists a backward bottleneck if and only if $r<M$ and we know $G$ is fragile if and only if $r<M$. Hence, $G$ is fragile if and only if there exists a backward bottleneck.
\end{proof}

\section*{Appendix B}\label{s:appendixB} 
In this section we explain how to find the kernel nodes of a large structured system and provide an example.
To compute the kernel nodes for a large structured system, first, compute its adjacency matrix. Then replace all the ones in the adjacency matrix with nonzero random numbers; denote the resulting matrix with $S$. Next compute a basis for the null space of $S$; in most computational languages, the output will be a matrix in which each column is a basis vector for the null space of the structure matrix; denote this matrix with $K$. By definition~\ref{defnn}(2), a regular node will be a node corresponding to a row that is zero (near machine zero) in every column, and a kernel node corresponds to a row that is nonzero (or far from zero) in at least one column. Assign a vector to matrix $K$, which includes the sum along the columns of $K$ and denote it by $I$. The nonzero entries of vector $I$ correspond to the kernel nodes since having a nonzero entry in $I$ means there has been a nonzero entry in the corresponding row in matrix $K$. Below we compute the set of kernel nodes for the structured system represented in Figure~\eqref{fg:4-troph d}(a) with python as an example. 

For the convenience of the reader, the structured graph in  Fig.~\ref{fg:4-troph}(a) is repeated below with its nodes labeled by numbers. First, 
the nonzero entries of the structure matrix is replaced by random numbers. Then using null\textunderscore space function from scipy module we compute a basis for the null space of the structured matrix $S$. Note that since the matrix is filled with random numbers, the output may vary in each run, but the output below shows what to expect. 
Each column of the output matrix below corresponds to a basis vector. 
By inspection of the output, it is noted that coordinates $1, 2, 7, 8, 9$, corresponding to the red nodes, have values at or near machine zero, in every column. For the convenience of the reader,  the rows of output that are near zero are colored red. Each row corresponding to a bottle node is far from zero at least in one of its entries. Hence, in the next step, we sum  along the columns of the matrix $K$. If a row is always near to zero, sum of entries of that row will remain very close to zero, otherwise, it will be a nonzero number far from zero. The entries far from zero correspond to kernel nodes which are labeled blue in the output below.
\begin{figure}[H]
\centering
\begin{tikzpicture}
    [ultra thick]
    \node[rdnode](a11)[xshift=-1cm,yshift=1cm]{1};
    \node[rdnode](a12)[xshift=+1cm,yshift=1cm]{2};
   \node[bdnode](a21)[xshift=-2.2cm,yshift=-0.5cm]{3};
    \node[bdnode](a22)[xshift=-.9cm,yshift=-0.5cm]{4};
    \node[bdnode](a23)[xshift=+.9cm,yshift=-0.5cm]{5};
    \node[bdnode](a24)[xshift=+2.2cm,yshift=-0.5cm]{6};
    \node[rdnode](a31)[xshift=-2cm,yshift=-2cm]{7};
    \node[rdnode](a32)[xshift=0cm,yshift=-2cm]{8};
    \node[rdnode](a33)[xshift=2cm,yshift=-2cm]{9};
   \node[bdnode](a41)[xshift=-3cm,yshift=-3.5cm]{10};
    \node[bdnode](a42)[xshift=-.95cm,yshift=-3.5cm]{11};
    \node[bdnode](a43)[xshift=+.95cm,yshift=-3.5cm]{12};
    \node[bdnode](a44)[xshift=+3cm,yshift=-3.5cm]{13};
\draw[<->,>=stealth,black] (a41)-- (a32);
\draw[<->,>=stealth,black] (a44)-- (a33);
\draw[<->,>=stealth,black] (a44)-- (a31);
\draw[<->,>=stealth,black] (a43)-- (a32);
\draw[<->,>=stealth,black] (a43)-- (a33);
\draw[<->,>=stealth,black] (a42)-- (a31);
\draw[<->,>=stealth,black] (a42)-- (a33);
\draw[<->,>=stealth,black] (a41)-- (a32);
\draw[<->,>=stealth,black] (a41)-- (a33);
\draw[<->,>=stealth,black] (a31)-- (a21);
\draw[<->,>=stealth,black] (a31)-- (a22);
\draw[<->,>=stealth,black] (a31)-- (a24);
\draw[<->,>=stealth,black] (a32)-- (a21);
\draw[<->,>=stealth,black] (a32)-- (a22);
\draw[<->,>=stealth,black] (a33)-- (a23);
\draw[<->,>=stealth,black] (a33)-- (a21);
\draw[<->,>=stealth,black] (a32)-- (a24);
\draw[<->,>=stealth,black] (a21)-- (a11);
\draw[<->,>=stealth,black] (a21)-- (a12);
\draw[<->,>=stealth,black] (a22)-- (a11);
\draw[<->,>=stealth,black] (a22)-- (a12);
\draw[<->,>=stealth,black] (a23)-- (a11);
\draw[<->,>=stealth,black] (a23)-- (a12);
\draw[<->,>=stealth,black] (a24)-- (a11);
\draw[<->,>=stealth,black] (a24)-- (a12);
\draw[<->,>=stealth,black](a31)to[bend right,out=110,in=100](a11);
\draw[<->,>=stealth,black] (a33)to[bend right,out=-100,in=-95](a12);
\node at ([xshift=.0cm,yshift=0.5cm]a11.north east){\Huge{\xmark}};
\node at ([xshift=.5cm,yshift=0.5cm]a11.north east){\Huge{\xmark}};
\node at ([xshift=1.cm,yshift=0.5cm]a11.north east){\Huge{\xmark}};
\end{tikzpicture}
\end{figure}
\newpage
\begin{lstlisting}[frame = single, basicstyle=\footnotesize, language=Python]
import numpy as np
from scipy.linalg import null_space
S = np.random.randn(13,13)
edges = [[0,0,1,1,1,1,1,0,0,0,0,0,0],\
         [0,0,1,1,1,1,0,0,1,0,0,0,0],\
         [1,1,0,0,0,0,1,1,1,0,0,0,0],\
         [1,1,0,0,0,0,1,1,0,0,0,0,0],\
         [1,1,0,0,0,0,0,0,1,0,0,0,0],\
         [1,1,0,0,0,0,1,1,0,0,0,0,0],\
         [1,0,1,1,0,1,0,0,0,0,1,0,1],\
         [0,0,1,1,0,1,0,0,0,1,0,1,0],\
         [0,1,1,0,1,0,0,0,0,1,1,1,1],\
         [0,0,0,0,0,0,0,1,1,0,0,0,0],\
         [0,0,0,0,0,0,1,0,1,0,0,0,0],\
         [0,0,0,0,0,0,0,1,1,0,0,0,0],\
         [0,0,0,0,0,0,1,0,1,0,0,0,0]]
S=S*edges
K= null_space(S)
K
\end{lstlisting}

}
\begin{Verbatim}[commandchars=\\\{\}]
array([[ {\red 1.20717887e-16,  6.84428086e-17, -3.11732981e-16}],
       [{\red -2.27911323e-16, -9.70940348e-17,  5.54609082e-16}],
      \blue [ -9.21205896e-02, -1.55525688e-01,  3.89155479e-01],
      \blue [  3.50849655e-01,  1.91779975e-01, -6.49053981e-01],
     \blue  [ -2.77855603e-01,  2.05519591e-01, -2.29308997e-01],
      \blue [  1.24862901e-02, -1.71947247e-01,  3.48715679e-01],
       [{\red -6.49478401e-17, -1.22157121e-17,  1.15851999e-17}],
       [{\red -4.43549872e-17, -5.08385029e-17,  1.65382089e-16}],
       [ {\red 1.01923948e-16,  1.80608967e-17, -1.10077705e-16}],
     \blue  [ -3.41302851e-01,  4.41425689e-01, -1.10191730e-01],
     \blue  [ -6.81666946e-01, -2.93761792e-01, -1.16448517e-01],
    \blue   [  2.34523019e-01, -7.24053614e-01, -2.94561356e-01],
     \blue  [  3.93573583e-01,  2.48639855e-01,  3.75010446e-01]]
       \black)
\end{Verbatim}
\newpage
\begin{lstlisting}[frame = single, basicstyle=\footnotesize, language=Python]
import numpy as np

ar = np.array(abs(K))

binar = ar > pow(10, -12)

int_ar = binar.astype(int)

I = np.sum(int_ar, axis=1)

I
\end{lstlisting}

\begin{Verbatim}[commandchars=\\\{\}]
array([{\red 0, 0,{\blue 3, 3, 3, 3,} 0, 0, 0,{\blue 3, 3, 3, 3}}])
\end{Verbatim}

\end{document}